\documentclass[11pt]{amsart}
\usepackage{amsmath}
\usepackage{amssymb,latexsym}
\usepackage{amscd}
\usepackage{xspace}
\usepackage{verbatim}
\usepackage[english,francais]{babel}
\begin{document}
\numberwithin{equation}{section}
\selectlanguage{english}
\title[Open and finely open sets]{Maximal Solutions for $-\Gd u+u^q=0$ in open and finely open sets}
\author{Moshe Marcus }
\address{Department of Mathematics, Technion
 Haifa, ISRAEL}
 \email{marcusm@math.technion.ac.il}
\author{ Laurent Veron}
\address{Laboratoire de Math\'ematiques, Facult\'e des Sciences
Parc de Grandmont, 37200 Tours, FRANCE}
\email{veronl@lmpt.univ-tours.fr}

\thanks{Both authors were partially sponsored by the French -- Israeli cooperation program through grant
No. 3-4299.
 The first author (MM) also
wishes to acknowledge the support of the Israeli Science Foundation
through grant No. 145-05. }
\subjclass[2000]{Primary:     Secondary: }
\date{\today}
\newcommand{\txt}[1]{\;\text{ #1 }\;}
\newcommand{\tbf}{\textbf}
\newcommand{\tit}{\textit}
\newcommand{\tsc}{\textsc}
\newcommand{\trm}{\textrm}
\newcommand{\mbf}{\mathbf}
\newcommand{\mrm}{\mathrm}
\newcommand{\bsym}{\boldsymbol}
\newcommand{\scs}{\scriptstyle}
\newcommand{\sss}{\scriptscriptstyle}
\newcommand{\txts}{\textstyle}
\newcommand{\dsps}{\displaystyle}
\newcommand{\fnz}{\footnotesize}
\newcommand{\scz}{\scriptsize}
\newcommand{\be}{\begin{equation}}
\newcommand{\bel}[1]{\begin{equation}\label{#1}}
\newcommand{\ee}{\end{equation}}
\newtheorem{subn}{\name}
\renewcommand{\thesubn}{}
\newcommand{\bsn}[1]{\def\name{#1$\!\!$}\begin{subn}}
\newcommand{\esn}{\end{subn}}
\newtheorem{sub}{\name}[section]
\newcommand{\dn}[1]{\def\name{#1}}   
\newcommand{\bs}{\begin{sub}}
\newcommand{\es}{\end{sub}}
\newcommand{\bsl}[1]{\begin{sub}\label{#1}}
\newcommand{\bth}[1]{\def\name{Theorem}\begin{sub}\label{t:#1}}
\newcommand{\blemma}[1]{\def\name{Lemma}\begin{sub}\label{l:#1}}
\newcommand{\bcor}[1]{\def\name{Corollary}\begin{sub}\label{c:#1}}
\newcommand{\bdef}[1]{\def\name{Definition}\begin{sub}\label{d:#1}}
\newcommand{\bprop}[1]{\def\name{Proposition}\begin{sub}\label{p:#1}}
\newcommand{\bnote}[1]{\def\name{\mdseries\scshape Notation}\begin{sub}\label{n:#1}}
\newcommand{\bproof}{\begin{proof}}
\newcommand{\eproof}{\end{proof}}
\newcommand{\bcom}{}
\newcommand{\req}{\eqref}
\newcommand{\rth}[1]{Theorem~\ref{t:#1}}
\newcommand{\rlemma}[1]{Lemma~\ref{l:#1}}
\newcommand{\rcor}[1]{Corollary~\ref{c:#1}}
\newcommand{\rdef}[1]{Definition~\ref{d:#1}}
\newcommand{\rprop}[1]{Proposition~\ref{p:#1}}
\newcommand{\rnote}[1]{Notation~\ref{n:#1}}
\newcommand{\BA}{\begin{array}}
\newcommand{\EA}{\end{array}}
\newcommand{\BAN}{\renewcommand{\arraystretch}{1.2}
\setlength{\arraycolsep}{2pt}\begin{array}}
\newcommand{\BAV}[2]{\renewcommand{\arraystretch}{#1}
\setlength{\arraycolsep}{#2}\begin{array}}
\newcommand{\BSA}{\begin{subarray}}
\newcommand{\ESA}{\end{subarray}}
\newcommand{\BAL}{\begin{aligned}}
\newcommand{\EAL}{\end{aligned}}
\newcommand{\BALG}{\begin{alignat}}
\newcommand{\EALG}{\end{alignat}}
\newcommand{\BALGN}{\begin{alignat*}}
\newcommand{\EALGN}{\end{alignat*}}
\newcommand{\note}[1]{\noindent\textit{#1.}\hspace{2mm}}
\newcommand{\Remark}{\note{Remark}}
\newcommand{\forevery}{\quad \forall}
\newcommand{\1}{\\[1mm]}
\newcommand{\2}{\\[2mm]}
\newcommand{\3}{\\[3mm]}
\newcommand{\set}[1]{\{#1\}}
\def\({{\rm (}}
\def\){{\rm )}}
\newcommand{\st}[1]{{\rm (#1)}}
\newcommand{\lra}{\longrightarrow}
\newcommand{\lla}{\longleftarrow}
\newcommand{\llra}{\longleftrightarrow}
\newcommand{\Lra}{\Longrightarrow}
\newcommand{\Lla}{\Longleftarrow}
\newcommand{\Llra}{\Longleftrightarrow}
\newcommand{\warrow}{\rightharpoonup}
\def\dar{\downarrow}
\def\uar{\uparrow}
\newcommand{\paran}[1]{\left (#1 \right )}
\newcommand{\sqrbr}[1]{\left [#1 \right ]}
\newcommand{\curlybr}[1]{\left \{#1 \right \}}
\newcommand{\absol}[1]{\left |#1\right |}
\newcommand{\norm}[1]{\left \|#1\right \|}
\newcommand{\angbr}[1]{\left< #1\right>}
\newcommand{\paranb}[1]{\big (#1 \big )}
\newcommand{\sqrbrb}[1]{\big [#1 \big ]}
\newcommand{\curlybrb}[1]{\big \{#1 \big \}}
\newcommand{\absolb}[1]{\big |#1\big |}
\newcommand{\normb}[1]{\big \|#1\big \|}
\newcommand{\angbrb}[1]{\big\langle #1 \big \rangle}
\newcommand{\thkl}{\rule[-.5mm]{.3mm}{3mm}}
\newcommand{\thknorm}[1]{\thkl #1 \thkl\,}
\newcommand{\trinorm}[1]{|\!|\!| #1 |\!|\!|\,}
\newcommand{\vstrut}[1]{\rule{0mm}{#1}}
\newcommand{\rec}[1]{\frac{1}{#1}}
\newcommand{\opname}[1]{\mathrm{#1}\,}
\newcommand{\supp}{\opname{supp}}
\newcommand{\dist}{\opname{dist}}
\newcommand{\sign}{\opname{sign}}
\newcommand{\diam}{\opname{diam}}
\newcommand{\q}{\quad}
\newcommand{\qq}{\qquad}
\newcommand{\hsp}[1]{\hspace{#1mm}}
\newcommand{\vsp}[1]{\vspace{#1mm}}
\newcommand{\prt}{\partial}
\newcommand{\sms}{\setminus}
\newcommand{\ems}{\emptyset}
\newcommand{\ti}{\times}
\newcommand{\pr}{^\prime}
\newcommand{\ppr}{^{\prime\prime}}
\newcommand{\tl}{\tilde}
\newcommand{\wtl}{\widetilde}
\newcommand{\sbs}{\subset}
\newcommand{\Sbs}{\Subset}
\newcommand{\sbeq}{\subseteq}
\newcommand{\ovl}{\overline}
\newcommand{\unl}{\underline}
\newcommand{\nin}{\not\in}
\newcommand{\pfrac}[2]{\genfrac{(}{)}{}{}{#1}{#2}}
\newcommand{\tin}{\to\infty}
\newcommand{\ind}[1]{_{\scriptscriptstyle #1}}
\newcommand{\chr}[1]{\chi\ind{#1}}
\newcommand{\rest}[1]{\big |\ind{#1}}
\newcommand{\Sol}[2]{\mathrm{Sol}\ind{#2}{#1}}
\newcommand{\wkc}{weak convergence\xspace}
\newcommand{\wrto}{with respect to\xspace}
\newcommand{\cons}{consequence\xspace}
\newcommand{\consy}{consequently\xspace}
\newcommand{\Consy}{Consequently\xspace}
\newcommand{\Essy}{Essentially\xspace}
\newcommand{\essy}{essentially\xspace}
\newcommand{\mnz}{minimizer\xspace}
\newcommand{\sth}{such that\xspace}
\newcommand{\ngh}{neighborhood\xspace}
\newcommand{\nghs}{neighborhoods\xspace}
\newcommand{\seq}{sequence\xspace}
\newcommand{\seqs}{sequences\xspace}
\newcommand{\sseq}{subsequence\xspace}
\newcommand{\ifif}{if and only if\xspace}
\newcommand{\suff}{sufficiently\xspace}
\newcommand{\abc}{absolutely continuous\xspace}
\newcommand{\sol}{solution\xspace}
\newcommand{\subss}{sub-solutions\xspace}
\newcommand{\subs}{sub-solution\xspace}
\newcommand{\supers}{super-solution\xspace}
\newcommand{\superss}{super-solutions  \xspace}
\newcommand{\Wlg}{Without loss of generality\xspace}
\newcommand{\wlg}{without loss of generality\xspace}
\newcommand{\locun}{locally uniformly\xspace}
\newcommand{\bvp}{boundary value problem\xspace}
\newcommand{\bvps}{boundary value problems\xspace}
\newcommand{\bdw}{\partial\Gw}
\newcommand{\Capq}{C_{2/q,q'}}
\newcommand{\Wq}{W^{2/q,q'}}
\newcommand{\Wqdual}{W^{-2/q,q}}
\newcommand{\Wqdb}{W^{-2/q,q}_+(\bdw)}
\newcommand{\sbsq}{\overset{q}{\sbs}}
\newcommand{\app}[1]{\underset{#1}{\approx}}
\newcommand{\suppq}{\mathrm{supp}^q_{\bdw}\,}
\newcommand{\convq}{\overset{q}{\to}}
\newcommand{\barq}[1]{\bar{#1}^{^q}}
\newcommand{\prtq}{\partial_q}
\newcommand{\tr}{\mathrm{tr}\,}
\newcommand{\Tr}{\mathrm{Tr}\,}
\newcommand{\trin}[1]{\mathrm{tr}\ind{#1}}
\newcommand{\Ltrin}[1]{\text{\rm L-tr}\ind{#1}}
\newcommand{\Mtrin}[1]{\text{\rm M-tr}\ind{#1}}
\newcommand{\qcl}{$q$-closed\xspace}
\newcommand{\qop}{$q$-open\xspace}
\newcommand{\gsmod}{$\gs$-moderate\xspace}
\newcommand{\gsreg}{$\gs$-regular\xspace}
\newcommand{\qreg}{quasi regular\xspace}
\newcommand{\qeq}{$q$-equivalent\xspace}
\newcommand{\ppf}{\underset{f}{\prec\prec}}
\newcommand{\ofrown}{\overset{\frown}}
\newcommand{\modcon}{\underset{mod}{\lra}}
\newcommand{\ugb}[1]{u\chr{\Gs_\gb(#1)}}
\def\RN{\BBR^N}
\def\loc{\ind{\rm loc}}
\def\qloc{\ind{\ell(2,q')}}
\def\bcom{}
\def\ga{\alpha}     \def\gb{\beta}       \def\gg{\gamma}
\def\gc{\chi}       \def\gd{\delta}      \def\ge{\epsilon}
\def\gth{\theta}                         \def\vge{\varepsilon}
\def\gf{\phi}       \def\vgf{\varphi}    \def\gh{\eta}
\def\gi{\iota}      \def\gk{\kappa}      \def\gl{\lambda}
\def\gm{\mu}        \def\gn{\nu}         \def\gp{\pi}
\def\vgp{\varpi}    \def\gr{\rho}        \def\vgr{\varrho}
\def\gs{\sigma}     \def\vgs{\varsigma}  \def\gt{\tau}
\def\gu{\upsilon}   \def\gv{\vartheta}   \def\gw{\omega}
\def\gx{\xi}        \def\gy{\psi}        \def\gz{\zeta}
\def\Gg{\Gamma}     \def\Gd{\Delta}      \def\Gf{\Phi}
\def\Gth{\Theta}
\def\Gl{\Lambda}    \def\Gs{\Sigma}      \def\Gp{\Pi}
\def\Gw{\Omega}     \def\Gx{\Xi}         \def\Gy{\Psi}

\def\CS{{\mathcal S}}   \def\CM{{\mathcal M}}   \def\CN{{\mathcal N}}
\def\CR{{\mathcal R}}   \def\CO{{\mathcal O}}   \def\CP{{\mathcal P}}
\def\CA{{\mathcal A}}   \def\CB{{\mathcal B}}   \def\CC{{\mathcal C}}
\def\CD{{\mathcal D}}   \def\CE{{\mathcal E}}   \def\CF{{\mathcal F}}
\def\CG{{\mathcal G}}   \def\CH{{\mathcal H}}   \def\CI{{\mathcal I}}
\def\CJ{{\mathcal J}}   \def\CK{{\mathcal K}}   \def\CL{{\mathcal L}}
\def\CT{{\mathcal T}}   \def\CU{{\mathcal U}}   \def\CV{{\mathcal V}}
\def\CZ{{\mathcal Z}}   \def\CX{{\mathcal X}}   \def\CY{{\mathcal Y}}
\def\CW{{\mathcal W}}
\def\BBA {\mathbb A}   \def\BBb {\mathbb B}    \def\BBC {\mathbb C}
\def\BBD {\mathbb D}   \def\BBE {\mathbb E}    \def\BBF {\mathbb F}
\def\BBG {\mathbb G}   \def\BBH {\mathbb H}    \def\BBI {\mathbb I}
\def\BBJ {\mathbb J}   \def\BBK {\mathbb K}    \def\BBL {\mathbb L}
\def\BBM {\mathbb M}   \def\BBN {\mathbb N}    \def\BBO {\mathbb O}
\def\BBP {\mathbb P}   \def\BBR {\mathbb R}    \def\BBS {\mathbb S}
\def\BBT {\mathbb T}   \def\BBU {\mathbb U}    \def\BBV {\mathbb V}
\def\BBW {\mathbb W}   \def\BBX {\mathbb X}    \def\BBY {\mathbb Y}
\def\BBZ {\mathbb Z}

\def\GTA {\mathfrak A}   \def\GTB {\mathfrak B}    \def\GTC {\mathfrak C}
\def\GTD {\mathfrak D}   \def\GTE {\mathfrak E}    \def\GTF {\mathfrak F}
\def\GTG {\mathfrak G}   \def\GTH {\mathfrak H}    \def\GTI {\mathfrak I}
\def\GTJ {\mathfrak J}   \def\GTK {\mathfrak K}    \def\GTL {\mathfrak L}
\def\GTM {\mathfrak M}   \def\GTN {\mathfrak N}    \def\GTO {\mathfrak O}
\def\GTP {\mathfrak P}   \def\GTR {\mathfrak R}    \def\GTS {\mathfrak S}
\def\GTT {\mathfrak T}   \def\GTU {\mathfrak U}    \def\GTV {\mathfrak V}
\def\GTW {\mathfrak W}   \def\GTX {\mathfrak X}    \def\GTY {\mathfrak Y}
\def\GTZ {\mathfrak Z}   \def\GTQ {\mathfrak Q}
\font\Sym= msam10
\def\SYM#1{\hbox{\Sym #1}}
\def\bmn{\mathbf{n}}
\def\bma{\mathbf{a}}
\def\prtn{\prt_{\bmn}}

\def\txin{\txt{in}}
\def\txon{\txt{on}}
\def\C2q{C_{2,q'}}
\def\W2q{W^{2,q'}}
\def\w02q{W_{0,\infty}^{2,q'}}
\def\qfine{\text{$\C2q$-finely}\xspace}
\def\qsupp{\mathrm{supp}\ind{(2,q')}}
\def\qcomp{\text{$\C2q$-pseudo compact}\xspace}
\def\qae{\text{$\C2q$ a.e.}\xspace}
\def\qinterior{\mathrm{int}_q}
\def\limsim{\overset{\scriptscriptstyle lim}{\sim}}
\def\limsbs{\overset{\scriptscriptstyle lim}{\subset}}
\def\qsim{\overset{\scriptscriptstyle \, q}{\sim}}
\def\qsbs{\overset{\scriptscriptstyle \, q}{\sbs}}
\def\zta{\gz_{\eta}}
\def\mq{2m/(q-1)}
\def\qtop{$\C2q$-fine topology\xspace}
\def\BPq{(B-P)$_q$\xspace}
\def\prtql{$\prt_q$-large\xspace}
\def\quasi{$\C2q$-quasi\xspace}
\def\qlim{\displaystyle \lim^q}
\def\qleq{\overset{\scriptscriptstyle \, q}{\leq}}

\begin{abstract}
\selectlanguage{english}
We derive sharp  estimates for the maximal solution
$U$ of (*) $-\Delta u+u^q=0$ in an  arbitrary
open set $D\subset\mathbb R^N$ The estimates
 involve the Bessel capacity $C_{2,q'}$, for $q$ in the
 supercritical range $q\geq q_{c}:=N/(N-2)$.
 We provide a  pointwise necessary and sufficient
condition, via a Wiener type criterion, in order that
$U(x)\to\infty$ as $x\to y$ for given $y\in\prt D$.
This completes the study of such  criterions carried out in
\cite{DL} and \cite{La}.
Further, we extend the notion of solution to $\C2q$ finely open sets
and show that, under very general conditions, a boundary value problem
with blow-up on a specific subset of the boundary is well-posed.
This implies, in particular, uniqueness of large solutions.

\vskip 0.5\baselineskip

\selectlanguage{francais}
\noindent{\bf Solutions maximales de $\Delta u=u^q$ dans des ensembles ouverts et finement ouverts}\\
{\sc R\'esum\'e.}\hskip 2mm
Nous d\'emontrons des estimations pr\'ecises pour la solution maximale $U$ de (*) $-\Delta u+u^q=0$ dans un domaine arbitraire $D\subset \BBR^N$. Ces estimations impliquent la capacit\'e de Bessel $C_{2,q'}$, pour $q$ appartenant \`a l'intervalle sur-critique $q\geq q_{c}:=N/(N-2)$. Nous donnons une condition n\'ecessaire et suffisante ponctuelle, via un crit\`ere de type Wiener, pour que $U(x)\to\infty$ quand $x\to y$ pour un $y\in\prt D$ arbitraire. Ce r\'esultat compl\`ete l'\'etude de tels crit\`eres men\'ee dans \cite{DL} et \cite{La}. En outre, nous \'etendons la notion de solution \`a des ensembles  finement ouverts pour la topologie $\C2q$ et montrons que, sous des conditions tr\`es g\'en\'erales, un probl\`eme aux limites avec explosion sur un sous-ensemble sp\'ecifique du bord est bien pos\'e. Cela implique en particulier l'unicit\'e des grandes solutions.
\end{abstract}
\maketitle
\tableofcontents
\section{Introduction}
\label{Introduction} \setcounter{equation}{0}
 In this paper we
study solutions of the equation
\begin{equation}\label{eq}
 -\Gd u+|u|^{q-1}u=0,
\end{equation}
in $\Gw\setminus F$, $\Gw$ a smooth domain in $\RN$, $N\geq 3$ and $F\sbs \Gw $, $F$  compact
or, more generally, a bounded set, closed in the $\C2q$ fine
topology. Here $q>1$ and $\C2q$ refers to the Bessel capacity with the specified indexes.
If $1<q<q_c=N/(N-2)$ then
the fine topology is equivalent to the Euclidean topology. Therefore, throughout the paper we shall assume that
$q\geq q_c$, in which case the two topologies are different.

If $D$ is an open set and $\mu$ is a Radon measure in $D$, a function
$u\in L^q\loc(D)$ is a solution of
\begin{equation}\label{eqmu} -\Gd
u+|u|^{q-1}u=\mu \txin{D}
\end{equation}
if the equation is satisfied in the distribution sense. It is known
\cite{BP} that \req{eqmu} possesses a solution \ifif $\mu$ vanishes
on sets of $\C2q$ capacity zero. When this is the case we say that
$\mu$ satisfies the \BPq condition (i.e., the Baras-Pierre
condition). If $D=\RN$ and $\mu$ is a Radon measure satisfying this
condition then \req{eqmu} possesses a {\em unique} solution.

Further, if $D$ is open, it is known that $C_{2,q'}(\RN\sms D)=0$ \ifif  the only solution of \req{eq} in $D$
is the trivial solution.
In view of the Keller -- Osserman estimates, the set of solutions of \req{eq} in $D$ (denoted by $\CU_D$) is
uniformly bounded in compact subsets of $D$
and every \seq of solutions possesses a \sseq which converges to a solution $u$.
Finally the compactness together with the maximum principle imply that $\max\,\CU_D$ is a solution in $D$.
The {\em maximal solution} in $D$ is denoted by $U_F$, $F=\RN\sms D$.

Now suppose that $F=\cup_{n=1}^\infty K_n$ where $\set{K_n}$ is an increasing \seq of
compact sets \sth $$\C2q(F\sms  K_j)\to 0.$$
 Then $\set{U_{K_n}}$ is non-decreasing
and we denote $V_F:=\lim U_{K_n}$. In this case $F$ may not be closed; in fact, it may be dense in $D=F^c$, so that
in general we cannot apply the
Keller -- Osserman estimates. Therefore, on this basis, it is not even clear whether
$V_F$ is finite a.e. in $D$. It will be shown in the course of this paper that this is actually the case.

Naturally, further questions come up:
Is $V_F$, in some sense, a generalized solution of \req{eq} in $D$ and, if so,
is it the maximal solution? Is it possible to characterize $V_F$ in terms of its behavior at the boundary?


The main objective of this paper is the study of properties of the
maximal solution of \req{eq} in $F^c$, {\em first} in the case that
$F$ is compact; {\em secondly} in the case that $F$ is merely \qfine
closed. In the second case we introduce a new notion of solution
which we call a {\em $\C2q$-strong solution} (see
\rdef{finesolution}) and show that  $V_F$ is indeed a solution in
this sense and that it is the maximal solution. We also show that
many of the properties of the set of classical solutions are shared
by the class of $\C2q$-strong solutions.

\bcom
If
$F$ is merely $\C2q$ finely closed, a new definition is needed,
because in this case $D$ need not be open. In fact $F$ may even be dense in $D$ (relative to Euclidean topology).


The following results are well known
 in the case that $F$ is {\em compact} (see \cite{BP}, \cite{LN}, \cite{MV1'} respectively). We denote by
 $\CU_D$ the set of non-negative solutions of \req{eq} in $D$.\2
{\bf (A) (removability)}   {\em If $F$ is a compact set \sth $\C2q(F)=0$ then the only  solution of \req{eq} in $D$
is the trivial solution.}

{\bf (B)(maximal solution)} {\em Let $F$ be a compact set. Then there exists a  solution $U_F$
of \req{eq} in $D$ \sth, for every $u\in \CU_D$, $U_F\geq u$ a.e. in $D$. In fact $U_F=\max\CU_D$.}

\note{Note} $U_F$ is called {\em the maximal solution} of \req{eq} in $D$.
{\bf (C) (compactness)} (i) {\em If $F$ is compact then $\CU_D$ is uniformly bounded in every compact subset of $D$.}\1
(ii) {\em If $\set{v_k}\sbs \CU_D$,
there exists a \sseq which converges pointwise to a solution $v\in \CU_D$.}

One of our objectives is to extend the notion of 'solution' to \qfine open sets and to show that these results, appropriately modified,
remain valid in this general case.

\begin{equation}\label{problem1}\BAL
-\Gd u+u^q&=0 &&\text{in}\q\Gw\sms F\\
u&=0 &&\text{on}\q\bdw \EAL
\end{equation}
when $\Gw$ is either a smooth bounded domain or $\Gw=\BBR^N$. This
solution will be denoted by $U_F^\Gw$ or simply by $U_F$, if
$\Gw=\BBR^N$. Note that $U_F^\Gw\leq U_F$.
\end{comment}
For $F$ {\em compact}, the properties of $U_F$ have been intensively investigated,
especially  in the last twenty years. A question that received special attention
 was the existence, uniqueness and estimates of
solutions of the \bvp
\begin{equation}\label{largeBVP}\BAL
   -\Gd u+|u|^{q-1}u&=0  \txin{D=F^c,} \\ \lim_{D\ni x\to y}u(x)&=\infty \forevery y\in \prt D.
\EAL\end{equation}
The question of existence reduces to the question whether $U_F$ blows up everywhere on the boundary.

\par A solution of \req{largeBVP} is called a {\em large solution} of \req{eq} in $D$. If $D$ is a smooth domain with compact boundary,
it is known that a large solution exists and is unique, (see \cite{LN}, \cite{BM1}, \cite {BM2}, \cite{Ve1}). These results were extended in various ways,
weakening
the assumptions on the domain, extending it to more general classes of equations and obtaining more information on the
asymptotic behavior of solutions at the boundary, (see  \cite{BM3}, \cite{MV0}, \cite{LM1}, \cite{BM4}  and references therein).

 In the present paper we also consider two
related notions:\1
(a) A solution $u$ is an {\em almost large solution}  of \req{eq} in $D$ if 
\begin{equation}\label{q-prtlarge}
 \lim_{D\ni x\to y}u(x)=\infty \q\text{\qae  $y\in \prt F$.}
\end{equation}
This notion is, in a sense, more natural, because (as we shall show) $U_F$ is invariable \wrto $\C2q$ equivalence of sets.
(Two Borel sets $E$ , $F$ are $\C2q$ equivalent if $\C2q(F\triangle E)=0$.)\1
 (b) A solution
$u$ of \req{eq} is a {\em $\prt_q$-large solution} in $D$ if
\begin{equation}\label{prtqlarge}
 \lim_{D\ni x\to y}u(x)=\infty \q\text{\qae  $y\in \prt_q F$,}
\end{equation}
where $\prt_q F$ denotes the boundary of $F$ in the \qtop.\1
\indent Here is a quick review of results pertaining to the case $F$ compact.

In the {\em subcritical case}, i.e. $1<q<q_c:=N/(N-2)$, the properties of $U_F$ are
well understood. In this case $C_{2,q'}(F)>0$ for any
non-empty set and it is classical that positive solutions may have
isolated point singularities of two types: weak and strong.
This easily implies that the maximal solution $U_F$ is always a
large solution in $F^c$. Sharp estimates of the large solution where obtained in \cite{MV4}. In addition it is proved in
\cite{Ve2} that
the large solution is unique if $\prt F^c\subset \prt \overline{F^c}^c$.


 In the subcritical case,  solutions with point singularities served as
  building blocks for solutions with general singularities. In the {\em supercritical case}, i.e. $q\geq q_{c}$, the situation is much
more complicated, because there are no solutions with point singularities.

 Sharp estimates for
 $U_F$ were obtained by Dhersin and Le Gall \cite{DL} in the case $q=2$, $N\geq 4$. These estimates
 were expressed in
terms of the Bessel capacity $C_{2,2}$
 and  were used to provide
 a Wiener type criterion -- to which we refer as (WDL; 2) -- for the pointwise blow up  of $U_F$, i.e., given $y\in F$,
\begin{equation}\label{pt blowup}
 \lim_{F^c\ni x\to y} U_F(x)= \infty \iff \text{the (WDL; 2) criterion is satisfied at y.}
\end{equation}
These results were obtained by probabilistic tools; hence the
restriction to $q=2$.

\par Labutin \cite{La}  extended the results of
\cite{DL} in the case $q> q_c$. Specifically, he obtained sharp estimates for $U_F$ similar to those in \cite{DL},
with $C_{2,2}$ replaced by $C_{2,q'}$.
These estimates were used to obtain a Wiener criterion
involving $\C2q$ (we refer to it as (WDL;q)) relative to which the following was proved:
\begin{equation}\label{blowup}
 \text{$U_F$ is a large solution } \iff  \text{ (WDL;q) holds {\em everywhere} in $F$.}
\end{equation}
\par Of course this result is weaker then \req{pt blowup}. However a
careful examination of Labutin's proof reveals that, in the case $q>q_c$, his argument
actually proves \req{pt blowup}. In the case $q=q_c$
Labutin's estimate was not sharp and it did not yield \req{pt
blowup} although it was sufficient in order to obtain \req{blowup}.

Uniqueness was not discussed in the above papers. Necessary and sufficient conditions are not yet known.
Sufficient conditions for uniqueness of  large solutions, for arbitrary $q>1$,
can be found in \cite{MV0}, \cite{MV2} and references therein. Uniqueness will also be one of the main subjects of
the present work.

The first part  of the present paper (Sections 2-4) is devoted to the study of the maximal solution $U_F$ when
{\em $F$ is  compact} and of the
almost large solution in bounded open sets.
Here is the list of main results obtained in this part of the paper:\2
\noindent {\bf I.} Sharp capacitary estimates of $U_F$ in the full
supercritical range $q \geq q_c$, $N\geq 3$. As a result, we show
that a variant of \req{pt blowup} holds in the
entire supercritical range.  Specifically, we show that, for $y\in F$,
\begin{equation}\label{q-pt blowup}
 \lim_{F^c\ni x\to y} U_F(x)= \infty \iff W_F(y)=\infty,
\end{equation}
where $W_F:\RN\to [0,\infty]$ is the  {\em capacitary potential of $F$,} (see \req{WF0} for its definition).

For $q>q_c$ the condition $W_F(y)=\infty$ is equivalent to the (WDL;q) criterion mentioned before.
However our proof does not require separate
treatment of the border case $q=q_c$ and is simpler than the proof in \cite{La} even for $q>q_c$. \1
{\bf II.} {\em For every compact set $F$, $U_F$ is an almost large solution in $F^c$ and $U_F$ is \gsmod.}
\par The statement '$U_F$ is \gsmod` means that there exists a monotone increasing \seq of bounded, positive measures
concentrated in $F$,
$\set{\mu_n}$, satisfying the \BPq condition, \sth $u_{\mu_n} \uar U_F$.

 Finally we establish an existence and uniqueness result; for its statement we need some additional notation. For
  any set $E\sbs \RN$
   $$\wtl E = \text{closure of $E$ in the $\C2q$-fine topology,} \q \prt_q E:=\wtl E\cap \wtl{E^c}.$$
 {\bf III.}  {\em Let $\Gw=\cup\Gw_n$, where
$\set{\Gw_n}$ is an increasing \seq of open sets, and put $D_n=\RN\sms \ovl\Gw_n$. Assume that
\begin{equation}\label{III-1}
 \C2q(\Gw\sms \Gw_n)\to 0 \txt{and}  \C2q(\prt\Gw_n\sms \wtl D_n)\to 0.
\end{equation}
Then the \bvp
\begin{equation}\label{bvp-aeinfty}
  -\Gd u+u^q=0 \txin{\Gw},\q \lim_{\Gw\ni x\to y}u(x)=\infty \q\text{for \qae  $y\in \prt_q\Gw$}
\end{equation}
possesses exactly one solution.\2
}
In other words, an open set $\Gw$ as above, possesses exactly one $\prt_q$-large solution. If $\prt \Gw$ is compact then,
this solution is an almost large solution. Indeed,
by {\bf II},
the maximal solution
$U_{\prt\Gw}$ is an almost large solution in $\Gw$.
Since $\prt_q\Gw\sbs \prt\Gw$, this implies that $U_{\prt\Gw}$ is a $\prt_q$-large solution. By  {\bf III},  $U_{\prt\Gw}$ is
the unique such solution in $\Gw$.\2
\indent In the second part of the paper (Sections 5-7) we extend our
investigation to the case where $F$ is $\C2q$ finely closed. We
introduce the notion of $\C2q$-strong solution in $D=\RN\sms F$,
which is now merely \qfine open, and prove that $V_F$ is a
$\C2q$-strong solution. By definition a $\C2q$-strong solution
belongs to a
 certain type of local Lebesgue space described in Section \ref{sec:weaksol}  below. Further we derive integral
 a-priori estimates which serve to replace the Keller-Osserman estimate in this case.
Using them we prove removability and compactness results. In
addition we show that the capacitary estimates  {\bf I} and the
Wiener criterion for pointwise blowup, namely \req{q-pt blowup},
persist for $V_F$. We also establish  the following version of {\bf
II}:\1 {\bf II'.} {\em For every \qfine closed set $F$, $V_F$ is the
maximal $\C2q$-strong solution in $F^c$. $V_F$ is a \prtql solution
and it is \gsmod.}\2
\indent Finally, we have the following  existence and uniqueness result:\1
{\bf III'.} {\em Let $\Gw$ be a \qfine open set. Let $\set{G_n}$ be a \seq of open sets \sth
\begin{equation}\label{III-2}
 \C2q(G_n\Gd \Gw)\to 0,\qq \C2q(\prt G_n\sms \prt_q\wtl G_n)\to 0.
\end{equation}
Then \req{bvp-aeinfty}  possesses exactly one $\C2q$-strong solution. The definition of blow up at the boundary
is defined in a manner appropriate for this class of solutions (see \rdef{boundary data})
\2}
Note that here we  do not assume that $G_n$ is contained in $\Gw$ or contains $\Gw$. If $\Gw\sbs G_n$ for every $n\in \BBN$ then
\req{III-2} implies \req{III-1}.

This seems to be the first study of the subject
in the setting of the $\C2q$ fine topology, introducing a notion of solution in sets where the classical distribution
derivative is not applicable. However the related subject of 'finely harmonic functions' has been studied
for a long time (see e.g. \cite{Fuglede}). Finely harmonic functions are defined on finely open sets relative to
classical $C_{1,2}$-capacity; however their definition depends on specific properties of harmonic functions
(e.g. the mean value property).

The framework presented here is particularly suitable for the study or \req{eq} and \req{eqmu} because
limits of solutions in open domains lead naturally to $\C2q$ strong solutions in \qfine open sets.
The underlying limit is relatively weak, namely, limit in the topology of a local Lebesgue space defined by a family
of weighted semi-norms
with weights in $W^{2,q'}(\RN)$ that are {\em bounded and compactly supported in the finely open set} (see Section \ref{sec:weaksol}).

At present this framework is presented mainly in the context of the study of maximal solutions and
uniqueness of solutions with blow up on the $\C2q$ boundary. A more detailed study, including an
 extension to more general boundary value problems will appear elsewhere.
 \vskip 3mm
\hspace{1cm}\underline{Partial list of notations}
\begin{itemize}
\item $[a<f<b]$ means $\set{x:\, a<f(x)<b}$.
\item $A\Gd B=(A\cup B)\sms (A\cap B)$.
\item If $f,g$ are non-negative functions with domain $D$ then $f\sim g$ means\\
 that there exits a constant $C$ \sth
$C^{-1}f\leq g\leq Cf$.
\item $A\qsim B$ means $\C2q(A\Gd B)=0$,\q $A\sbsq B$ means $\C2q(A\sms B)=0$.
\item $\wtl A$ means 'the closure of $A$ in the $\C2q$ fine topology'.
\item $\prt_q A$ means 'the boundary of $A$ in the $\C2q$ fine topology'.
\item $\qinterior A$ means 'the interior of $A$ in the $\C2q$ fine topology'.
\item $A\Sbs B$ means '$A$ bounded and $\ovl A\sbs B$\,'.
\item $B_r(x_0)=\set{x\in \RN:\,\absol{x-x_0}<r}$.
\item $\chr{A}$ denotes the characteristic function of the set $A$.
\item \BPq condition: A  measure $\mu$ satisfies this condition if $|\mu|(E)=0$ for every
Borel set $E$ \sth $\C2q(E)=0$.
\item $u_\mu$ denotes the solution of \req{eqmu} in $\RN$ when $\mu$ is a Radon\\
measure satisfying the \BPq condition.
\end{itemize}

\section{Upper estimate of the maximal solution.}
In this section $F$ denotes a non-empty compact set in $\RN$ and the
maximal solution of \req{eq} in $\RN\sms F$ is denoted by $U_F$.
Further, for $x\in \RN$, we denote
\begin{equation}\label{Fm} \BAL &T_m(x)=\set{y\in \BBR^N:2^{-(m+1)}\leq \absol{y-x}\leq 2^{-m}},\\
  &F_m(x)=F\cap T_m(x), \q F^*_m(x)=F\cap \ovl B_{2^{-m}}(x),\EAL
\end{equation}
\begin{equation}\label{WF0}\BAL
 W_F(x)&= \sum_{-\infty}^{\infty} 2^{\frac{2m}{q-1}}C_{2,q'}\left(2^{m}F_m(x)\right),\\
 W^*_F(x)&= \sum_{-\infty}^{\infty} 2^{\frac{2m}{q-1}}C_{2,q'}\left(2^{m}F^*_m(x)\right).\EAL
\end{equation}
We call $W_F$ the {\em $C_{2,q'}$-capacitary potential} of $F$. It is known that the  two functions in \req{WF0} are equivalent,
i.e., there exists a constant $C$ depending only on $q,N$ such that
\begin{equation}\label{WF1}
 W_F(x)\leq W^*_F(x)\leq CW_F(x)
\end{equation}
see e.g.  \cite{MV3}.

If $K$ is a compact subset of a domain $\Gw$ put,
\begin{equation}\label{XK}
 X_K(\Gw):=\set{\eta\in C_c^2(\Gw): 0\leq \eta\leq 1,\; \eta=1 \text{ on } N^K_\eta},
\end{equation}
where $N^K_\eta$ denotes an open \ngh of $K$ depending on $\eta$.

The following theorem is due to Labutin \cite{La}:
\bth{up-est} Let $q\geq q_c$. There exists a constant $C$ depending only on $q,N$ \sth, for every compact set $F$,
\begin{equation}\label{UF<WF}
   U_F(x)\leq C W_F(x) \forevery x\in D.
\end{equation}
\es
For the convenience of the reader we provide a concise proof; components of this proof will also be used later on in the paper.
The main ingredient in the proof is contained in the lemma stated below.

\bcom
Let $\Gw$ be a bounded smooth domain and denote by $\vgf\ind{\Gw}$  the solution of
\begin{equation}\label{specialsol}
 -\Gd \vgf=\chr{\Gw} \q\text{in }\RN,\q \lim_{|x|\tin}\vgf(x)=0
\end{equation}
and by  $\vgf^0\ind{\Gw}$ the solution of
\begin{equation}\label{specialsol0}
 -\Gd \vgf=1 \q\text{in }\Gw,\q \vgf(x)=0 \q\text{on }\bdw.
\end{equation}
\end{comment}

\blemma{intUK} Let $R>1$  and denote by $\vgf\ind{R}$  the solution of
\begin{equation}\label{specialsol}
 -\Gd \vgf=\chr{B_R(0)} \q\text{in }\RN,\q \lim_{|x|\tin}\vgf(x)=0.
\end{equation}
Given $\eta\in W^{2,q'}(\RN)$, $0\leq \eta\leq 1$,  put
$$ \zta=\vgf\ind{R}(1-\eta)^{2q'}.$$
There exists  a constant $\bar c(N,q,R)$ \sth, for every compact set $K\sbs B_1(0)$,
\begin{align}\label{UK1}
 &\int_{\RN\sms K}U_K^q\zta\,dx \leq \bar c\norm{\eta}^{q'}\ind{W^{2,q'}(\RN)}&& \forall \eta\in X_K(\RN),\\
 &\int_{B_R(0)\sms K}U_K(1-\eta)^{2q'}dx\leq \bar c\norm{\eta}^{q'}\ind{W^{2,q'}(\RN)}&& \forall \eta\in X_K(\RN).\label{UK1'}
\end{align}
\es
\bproof For $|x|\geq R+2$,
\begin{equation}\label{Wiener5'}\BAL
0<&\vgf\ind{R}(x)+U_K(x)\leq c|x|^{2-N}, && |\nabla\vgf\ind{R}(x)|+|\nabla U_K(x)|\leq c|x|^{1-N}
\EAL \ee
where $c=c(N,q,R)$.
For every $R'>R$ and $\eta\in W^{2,q'}(\RN)$,
\begin{equation}\label{Wiener5}
 \int_{B_{R'}(0)\sms K}\big( -U_K\Gd \zta+ U_K^q\zta\big)dx=-\rec{R'}\int_{\prt B_{R'}}(U_K\nabla\zta-\zta\nabla U_K)\cdot x dS.
\end{equation}
By \req{Wiener5'}, the right hand side of \req{Wiener5} tends to zero as $R'\tin$ and we obtain,
\begin{equation}\label{Wiener5+}
 \int_{D}\big( -U_K\Gd \zta+ U_K^q\zta\big)dx=0,
 \ee
 where $D:=\RN\sms K$.
Further,
\[\Gd\zta=\vgf\ind{R}\Gd(1-\eta)^{2q'}-(1-\eta)^{2q'}\chr{B_2}+2\nabla\vgf\ind{R}\cdot\nabla (1-\eta)^{2q'}\]
so that,
\begin{equation}\label{Wiener6}\BAL
&\int_D U_K^q\zta\,dx+\int_{B_R(0)\sms K}U_K(1-\eta)^{2q'}dx= \\
&\int_D U_K\big(\vgf\ind{R}\Gd((1-\eta)^{2q'})+2\nabla\vgf\ind{R}\cdot\nabla ((1-\eta)^{2q'})\big)dx.
\EAL\end{equation}
Now,
\[ \Gd((1-\eta)^{2q'})=-2q'(1-\eta)^{2q'-1}\Gd\eta+ 2q'(2q'-1)(1-\eta)^{2q'-2}|\nabla\eta|^2,\]
so that
\begin{equation}\label{Wiener7}
 \int_D U_K \vgf\ind{R}\Gd((1-\eta)^{2q'})dx\leq c(I_1+I_2),
\end{equation}
where
\[I_1:=\int_D U_K \vgf\ind{R}(1-\eta)^{2q'-1}|\Gd\eta|dx,\q I_2:=\int_D U_K \vgf\ind{R}(1-\eta)^{2q'-2}|\nabla\eta|^2dx. \]
The estimate of $I_1$  is standard.
\begin{equation}\label{Wiener8}
\BAL I_1&\leq \Big(\int_D U_K^q\zta\,dx\Big)^{1/q}\Big(\int_D\vgf\ind{R}(1-\eta)|\Gd\eta|^{q'}dx\Big)^{1/q'}\\
&\leq c\Big(\int_D U_K^q\zta\, dx\Big)^{1/q}\norm{\eta}\ind{W^{2,q'}(\RN)}.
\EAL \ee
To estimate $I_2$ we consider $\eta\in X_K(B_R(0))$ and use the interpolation inequality
\begin{equation}\label{gradphi2}
 \norm{|\nabla \eta|^2}\ind{L^{q'}(D)}\leq c(q,N,R)\norm{\eta}\ind{L^{\infty}(D)}\norm{D^2\eta}\ind{L^{q'}(D)}.
\end{equation}
We obtain,
\begin{equation}\label{Wiener9}\BAL
I_2&\leq \Big(\int_D U_K^q\zta\, dx\Big)^{1/q}\Big(\int_D\vgf\ind{R} |\nabla\eta|^{2q'}dx\Big)^{1/q'}\\
&\leq c \Big(\int_D U_K^q\zta\, dx\Big)^{1/q}\norm{|\nabla\eta|^2}\ind{L^{q'}(D)}\\
&\leq c\Big(\int_D U_K^q\zta\, dx\Big)^{1/q}\norm{\eta}\ind{W^{2,q'}(\RN)}.
\EAL
\ee
for $\eta\in X_K(B_R(0))$.
Next
\begin{equation}\label{Wiener10}\BAL
&\int_D U_K\nabla\vgf\ind{R}\cdot\nabla ((1-\eta)^{2q'})dx
\leq 2q'\int_D U_K |\nabla\vgf\ind{R}||\nabla \eta|(1-\eta)^{2q'-1}dx\\
&\leq c\Big(\int_D U_K^q\zta\,dx\Big)^{1/q}\Big(\int_D \vgf\ind{R}^{-\frac{q'}{q}}(|\nabla\vgf\ind{R}||\nabla \eta|)^{q'}dx\Big)^{1/q'}.
\EAL
\end{equation}
In view of the fact that,  for $|x|\geq R+2$, $\vgf\ind{R} (x)\geq c|x|^{2-N}$, \req{Wiener5'} implies
$$\vgf\ind{R}^{-\frac{q'}{q}}|\nabla\vgf\ind{R}|^{q'}\leq c(N,q,R).$$
Hence
\begin{equation}\label{Wiener11}
\int_D U_K\nabla\vgf\ind{R}\cdot\nabla ((1-\eta)^{2q'})dx\leq c\Big(\int_D U_K^q\zta\,dx\Big)^{1/q}\norm{\eta}\ind{W^{1,q'}(\RN)}
\ee
Combining \req{Wiener6}--\req{Wiener11} we obtain \req{UK1} and \req{UK1'} for $\eta\in X_K(B_R(0))$.

Pick $\gw\in C_c^\infty(B_R(0)$ \sth
$0\leq \gw\leq 1$ and $\gw=1$ in $B_1(0)$. Given $\eta\in X_K(\RN)$, \req{UK1} and \req{UK1'} are valid if
$\eta$ is replaced by $\gw\eta$. However $(1-\eta)\leq (1-\gw\eta)$ and
$$\norm{\gw\eta}\ind{W^{2,q'}(\RN)}\leq c(N,q,\gw) \norm{\eta}\ind{W^{2,q'}(\RN)}.$$
Therefore \req{UK1} and \req{UK1'} are valid for every $\eta\in X_K(\RN)$.
\eproof
\bcor{intUK1} Assume that $R>3/2$.
There exists  a constant $c_1= c_1(N,q,R)$ \sth, for every
compact set $K\sbs  B_1(0)$
\begin{equation}\label{UKi1}
\int_{[3/2<|x|]}U_K^q\vgf\ind{R}\,dx+ \int_{{[3/2<|x|<R]}}U_K dx\leq c_1 C_{2,q'}(K)
\ee
\bcom
\begin{equation}\label{UKi2}
\int_{[3/2<|x|]}U_K (|(\Gd\eta|+|\nabla \eta|^2) \vgf\ind{R}\,dx\leq c_1 C_{2,q'}(K) \forevery \eta\in X_K(\RN)
\ee
\end{comment}
and
\begin{equation}\label{UK<CK}
\sup_{[3/2<|x|<R]}U_K\leq c_1 C_{2,q'}(K).
\ee
\es

\bproof
 Recall that
\begin{equation}\label{cap-def}
   C_{2,q'}(K)=\inf \set{\norm{\eta}^{q'}\ind{W^{2,q'}(\RN)}:\eta\in X_K(\RN)}.
\end{equation}
Let $\gw\in C_c^\infty(B_{3/2}(0))$ be a function \sth $0\leq \gw\leq 1$ and $\gw=1$ on $B_1(0)$.
For every compact set $K\sbs B_1(0)$ put
\begin{equation}\label{cap-def'}
   C_{2,q'}^\gw(K)=\inf \set{\norm{\gw\eta}^{q'}\ind{W^{2,q'}(\RN)}:\eta\in X_K(\RN)}.
\end{equation}
Clearly $C_{2,q'}(K)\leq  C_{2,q'}^\gw(K)$ and since
$$\norm{\gw\eta}^{q'}\ind{W^{2,q'}(\RN)}\leq c(N,q,\gw)\norm{\eta}^{q'}\ind{W^{2,q'}(\RN)}$$
we have
\begin{equation}\label{cap-equiv}
 C_{2,q'}(K)\leq  C_{2,q'}^\gw(K)\leq c(N,q,\gw)C_{2,q'}(K).
\end{equation}
Let $\set{\eta_n}$ be a \seq in $X_K(\RN)$ \sth
$$\norm{\gw\eta_n}^{q'}\ind{W^{2,q'}(\RN)}\to C^\gw_{2,q'}(K).$$
For $K\sbs B_1(0)$ \req{UK1} implies that,
\begin{equation}\label{UK3}\BAL
\int_{\RN\sms B_{3/2}}U_K^q\vgf\ind{R}\,dx&\leq \liminf_{n\tin}\int_{\RN\sms K}U_K^q\vgf\ind{R}(1-\gw\eta_n)^{2q'}\,dx\\
&\leq  c(N,q,\gw) C_{2,q'}(K).
\EAL\end{equation}
This proves \req{UKi1}. Inequality \req{UK<CK} (with the supremum over a slightly smaller annulus, say,
$[3/2+\ge<|x|<R-\ge]$ with $\ge>0$ \sth $R>3/2+2\ge$)
follows from \req{UKi1} and Harnack's inequality applied as in \cite{MV4}.
\eproof


 \noindent {\em Proof of Theorem \ref{t:up-est}.} Inequality
\req{UK<CK} implies,
\begin{equation}\label{UF<CF}
   U_F(x)\leq c(N,q)\gr\ind{F}(x)^{-2/(q-1)}\C2q\big(F/\gr\ind{F}(x)\big)
\end{equation}
for every  every compact set $F\sbs \RN$ and every $x\in \RN\sms F$ \sth $\gr_F(x)\geq (3/2)\diam F$.
Recall that $\gr\ind{F}(x):=\dist(x,F)$.

The implication relies on the similarity transformation associated with \req{eq}.
For any $a>0$, we have
\begin{equation}\label{similar}
U_F(x)=a^{-2/(q-1)}U_{F/a}(x/a)\forevery x\in \RN\sms F.
\end{equation}
Assume, as we may, that $F\sbs B_R(0)$, $R=\diam F$. Fix a point $\bar x\in \RN\sms F$ \sth $a:=\gr_F(\bar x)\geq R$.
Applying \req{UK<CK}
to the set $K=3F/2a$, we obtain
\[\BAL
U_F(\bar x)&= (2a/3)^{-2/(q-1)}U_{K}(3\bar x/2a)\\
&\leq c(N,q)a^{-2/(q-1)}\C2q(K)\leq c'(N,q)a^{-2/(q-1)}\C2q(F/a).
\EAL\]
\vskip 1mm

Next we show that \req{UF<CF} is equivalent to \req{UF<WF}.
Let $x\in D$ and put
\begin{equation}\label{M(x)}
    M(x):=min\set{m\in \BBN:2^{-m}<\gr\ind{F}(x)}.
\end{equation}
 Then $F_{k}(x)=\ems$ for all $k\geq M(x)$ and consequently
$$W_F(x)=\sum_{k=-\infty}^{M(x)} 2^{\frac{2k}{q-1}}C_{2,q'}\left(2^{k}F_k(x)\right)\leq
C 2^{\frac{2M(x)}{q-1}}\sup_{k\leq M(x)}\C2q(2^kF_k(x)).$$
However it is known that there exists a constant $A$ depending only on $q,N$ \sth
\begin{equation}\label{C(aE),1}
  \C2q(aE)\leq Aa^{N-\frac{2}{q-1}}\C2q(E) \forevery a\in (0,1),
\end{equation}
(see e.g. \cite{MV3}). In addition, for every $\ell>1$ there exists a constant $A$, depending on $q,N, \ell$, \sth
\begin{equation}\label{C(aE),2}
  \C2q(aE)\leq Aa^{N-\frac{2}{q-1}}\C2q(E) \forevery a\in (1,\ell).
\end{equation}
Inequality \req{C(aE),1} implies that
$$\BAL
W_F(x)&\leq C_1 2^{\frac{2M(x)}{q-1}}\C2q(2^{M(x)}F)\\
&\leq C_2\gr\ind{F}(x)^{\frac{-2}{q-1}}\C2q(2F/\gr\ind{F}(x))\leq
C_3\gr\ind{F}(x)^{\frac{-2}{q-1}}\C2q(F/\gr\ind{F}(x)),
\EAL$$
where $C_i$ are constants depending only on $q,N$. Thus \req{UF<WF} implies \req{UF<CF}.

To prove the implication in the opposite direction we use the following facts:\2
For every compact set $F$ there exists a \seq of bounded domains $\set{D_n}$ \sth
\begin{equation}\label{exhaust1}
 \text{(i) $\cup D_n=D:=F^c$, (ii) $\ovl D_n\sbs D_{n+1}$, (iii) $\prt D_n$ is Lipschitz.}
\end{equation}
Such a \seq is called a {\em Lipschitz exhaustion of} $D$.

If $u_n$ denotes the maximal solution of \req{eq} in $D_n$ then $u_n$ is the unique large solution of \req{eq} in $D_n$ (see \cite{MV2}),
$u_n> u_{n+1}$ in $D_n$ and $U_F=\lim u_n$.\2
\indent Let $E_i$, $i=1,\ldots,k$ be compact sets and $E:=\cup_1^kE_i$. One can choose a Lipschitz exhaustion
$\set{D_{i,n}}_{n=1}^\infty$ of
$D_i:=E_i^c$, $i=1,\ldots,k$, \sth the \seq $\set{D_n}$, $D_n=\cap_{i=1}^k D_{i,n}$, is a Lipschitz exhaustion of $D$. Let $u_{i,n}$
be the large solution in $D_{i,n}$. Then $v_n=\max(u_{1,n},\ldots,u_{k,n})$ is a subsolution while
$w_n=\sum_{i=1}^k u_{i,n}$ is a supersolution of \req{eq} in $D_n$. Hence  $u_n$, the unique large solution of \req{eq} in $D_n$,
satisfies $v_n\leq u_n\leq w_n$. \Consy
\begin{equation}\label{U<sumUi_}
  \max(U_{E_1},\ldots,U_{E_k})\leq U_E\leq \sum_{i=1}^k U_{E_i}.
\end{equation}

\bcom
Let $\set{E_i}$ be an increasing \seq of compact sets \sth $E:=\cup_1^\infty E_i$ is compact and
$K_i:=\overline{E\sms E_i}\downarrow \ems$. Then $U_{K_i}\dar 0$ and, by \req{U<sumUi_},
$$U_{E_i}\leq U_E\leq U_{E_i}+ U_{K_i}.$$
Hence
$$U_E=\lim U_{E_i}.$$
\end{comment}

Returning to the notation of \rth{up-est}, fix $\bar x\in D$ and put
$$i(\bar x)=\max\set{ i\in \BBZ: F\sbs \ovl B_{2^{-i}}(\bar x)}.$$
Then $F=\cup_{i(\bar x)}^{M(\bar x)} F_m(\bar x)$ and, by \req{U<sumUi_} and \req{UF<WF},
$$U_F \leq \sum_{m=i(\bar x)}^{M(\bar x)}U_{ F_i(\bar x)}\leq
C\sum_{m=i(\bar x)}^{M(\bar x)} 2^{\frac{2m}{q-1}}C_{2,q'}\left(2^{m}F_m(\bar x)\right).$$
In particular, $U_F(\bar x)\leq C W_F(\bar x)$. Thus \req{UF<CF} implies \req{UF<WF}.
\section {Lower estimate of the maximal solution}

We need the  following well-known result:
\bprop {BBP} Let $\gm$ be a positive measure in
$W^{-2,q}\ind{\loc}(\BBR^N)$ and
 let $\Gw$ be a smooth domain with compact boundary. Then there exists a unique solution
 of each of the  problems
\begin{equation}\label{equ-meas-0}
-\Gd u+u^q=\gm \txt {in }\Gw,\q u=0\txt{on} \bdw
\end{equation}
and
\begin{equation}\label{equ-meas-inf}
-\Gd u+u^q=\gm \txt {in }\Gw,\q u=\infty\txt{on} \bdw.
\end{equation}
If\, $\Gw$ is the whole space then there exists a unique solution $u_\mu$ of the equation
\begin{equation}\label{equ-meas}
-\Gd u+u^q=\gm\quad\text {in }\BBR^N.
\end{equation}
In each case the solution increases monotonically with $\mu$.
Finally $$u_\mu=\lim_{R\tin}u^R_{\mu,0}=\lim_{R\tin}u^R_{\mu,\infty}$$ where
$u^R_{\mu,0}$ and $u^R_{\mu,\infty}$ are the solutions  of \req{equ-meas-0} and \req{equ-meas-inf}
respectively, when $\Gw=B_R(0)$.
\es

When $\gm\in L^1_{loc}(\BBR^N)$  the result is due to Brezis \cite {Br1} and Brezis-Strauss \cite{BS}.
In the case of a smooth bounded domain $\Gw$, with $\gm\in W^{-2,q}(\Gw)$, the result is due to Baras
 and Pierre \cite {BP}.  The final observation is easily verified.

In this section, the solution of \req{equ-meas-0} will be denoted by $u\ind{\mu,\Gw}$.

If $F$ is a compact subset of $\BBR^N$, we define
\begin {eqnarray}\label {ca4}
V_F\,:=\sup\{u_{\gm}:\gm\in\mathfrak M_{+}(\BBR^N)\cap
W^{-2,q}(\BBR^N),\,\gm (F^c)=0\}.
\end {eqnarray}
Then $V_F$ is the maximal $\gs$-moderate solution of (\ref {eq}) in $F^c:=\BBR^N\setminus F$.
Obviously,
\begin{equation}\label{elementary1}
  V_F\leq U_F.
\end{equation}
\bcom
Assuming that $F\sbs B_1=B_1(0)$ and $R>1$ we denote
\begin{equation}\label{ca5}
 \unl U^R_{F}=\sup\{u^R_{\gm,0}:\gm\in\mathfrak M_{+}(\BBR^N)\cap
W^{-2,q}(\BBR^N),\,\gm (F^c)=0\}.
\end{equation}

\end{comment}

 We derive a lower estimate for $V_F$,  equivalent to the
upper estimate for $U_F$ obtained in the previous section. More precisely:
\bth{low-est} Assume that $F$ is a compact subset of $B_a(0)$ and let $D$ be a bounded smooth domain \sth $B_{6a}(0)\sbs D$.
 Then, for every $x\in B_{2a}(0)\sms F$,
there exists a positive measure $\gm^{x}\in W^{-2,q}(\RN)$,
supported in $F$, such that
\begin {equation}\label {WF<u}
cW_F(x)\leq u\ind{\gm^x,D}(x)\leq V_F\,(x),
\end {equation}
where $c$ is a positive constant depending only on $N,q$. In particular,
\begin{equation}\label{WF<lowU}
  c(N,q)W_F(x)\leq V_F\,(x) \forevery x\in \RN\sms F.
\end{equation}
\es
\bproof
Let  $\gl$ be a bounded Borel measure supported in $D$.
We denote by $\BBG_D[\gl]$ the Green potential of the measure in $D$:
\begin{equation}\label{R-Dir}
\BBG_D[\gl](\cdot):=\int_{D} g\ind{D}(\cdot,\gx)\,d\gl(\gx),
\end{equation}
where $g\ind{D}$ denotes Green's function in $D$.

If $\mu$ is a positive measure in $W^{-2,q}(D)$ then,
$$u\ind{\mu,D}\leq \BBG_D[\mu]$$
 and \consy
\begin{equation}\label{Gmu}
   u\ind{\mu,D}=\BBG_D[\mu]-\BBG_D[u\ind{\mu,D}^q]\geq \BBG_D[\mu]-\BBG_D[(\BBG_D[\mu])^q].
\end{equation}


Given $x_0\in B_{2a}\sms F$ we construct a measure
$\gm^{x_0}\in W^{-2,q}(\RN)$, concentrated on $F$ \sth \req{WF<u} holds.
By shifting the origin
to $x_0$ we may assume that $x_0=0$.
We observe that \req{WF<u} is invariant \wrto dilation. Therefore we may assume that $a=1/2$.
Following the shift and the dilation we have
\begin{equation}\label{ass1}
F\sbs B_{1}(0),\q B_{2}(0)\sbs D,\q 0\in F^c,
\end{equation}
and we have to prove \req{WF<u}, with an appropriate measure $\mu^0$, at $x=0$. The right inequality
in \req{WF<u} is trivial. Therefore we have to prove only that, for some non-negative measure
$\mu^0\in W^{-2,q}(\RN)$
supported in $F$,
\begin{equation}\label{WF<u;D}
 c(N,q)W_F(0)\leq u\ind{\mu^0,D}(0).
\end{equation}
In view of \req{ass1},  $$u\ind{\mu^0,B_2(0)}\leq u\ind{\mu^0,D}.$$
Therefore it is enough to prove \req{WF<u;D} for $D=B_2(0)$ which we assume in the rest of the proof.

In what follows we shall freely use the notation introduced in the previous section and write simply $F_n,\, T_n$  instead of
 $F_n(0), T_n(0)$ etc.\,. Observe that in the present case $F_n=\ems$ for $n\leq -1$ and $F_n^*=F$
 for $n\leq 0$.
For every non-negative integer $n$, let $\gn_{n}$ denote the capacitary measure of $2^nF_{n}$.
Thus, $\nu_n$ is a positive measure in $W^{-2,q}(\RN)$
 supported in $2^nF_n$ which satisfies
 \begin{equation}\label{nu1}
 \gn_{n}(2^nF_{n})=C_{2,q'}(2^nF_{n})=\norm {\gn_{n}}_{W^{-2,q}}^q.
\end{equation}
Let $\gm_{n},\gm$ be the Borel measures in $\RN$ given by
\begin{equation}\label{mu1}
\mu_{n}(A)=2^{-n(N-2q')}\nu_{n}(2^nA)  \q n=0,1,2,\ldots \q \mu=\sum_0^\infty \mu_n.
\end{equation}
Thus
\begin{align}\label{mu2}
 &\supp \mu_n\sbs F_n,&&\supp\mu\sbs F,\\
 &\mu_n(F_n)=2^{-n(N-2q')}C_{2,q'}(2^nF_{n}), && \mu\in W^{-2,q}(\RN).\label{mu3}
\end{align}
Observe also that, for $x,\gx\in B_1(0)$,
\begin{equation}\label{g-est}
   g\ind{D}(x,\gx)\approx\absol{x-\gx}^{2-N}.
\end{equation}
The notation $f\approx h$ means that there exists a positive constant $c$ depending only on $N,q$
\sth $c^{-1}h\leq f\leq ch$.

The remaining part of the proof consists of a series of estimates of the terms on the
right hand side of \req{Gmu} for  $\mu$ as above.\2
\note{Lower estimate of $\BBG_D[\mu]$ } Using \req{mu3} and \req{g-est} we obtain,
$$c_N2^{-(n+1)(2-N)}\leq g(0,\gx)\forevery \gx\in B_1(0),$$
\begin{equation}\label{G[mu]}\BAL
 \BBG_D[\mu](0)&=\sum_{n\geq 0}\int_{F_{n}}g(0,\gx)d\mu_{n}(\gx)\geq
c\sum_{n\geq 0}\int_{F_{n}}2^{n(N-2)}\,d\mu_{n}(\gx)\\
&=\sum_{n\geq 0}c2^{-2n/(q-1)}C_{2,q'}(2^nF_{n})=cW_F(0).
\EAL\end{equation}
\note{Upper estimate of $\BBG_D[(\BBG_D[\gm])^q](0)$} We prove that
\begin{equation}\label{upest0}\BAL
&  \BBG_D[(\BBG_D[\gm])^q](0)=\int_D g\ind{D}(0,\gx)\BBG_D[\mu]^q(\gx)d\gx\\
 & =\sum_{-1}^\infty\int_{T_k}g\ind{D}(0,\gx)\Big(\sum_{n\geq 0}\BBG_D[\mu_n](\gx)\Big)^q d\gx
 \leq c(N,q) W_F(0).
\EAL\end{equation}
This estimate requires several steps. Denote
\begin{align}\label{upest1}
I_1&=\sum_{k=3}^\infty\int_{T_k}g\ind{D}(0,\gx)\Big(\sum_{n=0}^{k-3}\BBG_D[\mu_n](\gx)\Big)^q\,d\gx\\
I_2&=\sum_{-1}^\infty\int_{T_k}g\ind{D}(0,\gx)
\Big(\sum_{n>k+2}\BBG_D[\mu_n](\gx)\Big)^q\,d\gx \label{upest2}\\
I_3&=\sum_{-1}^\infty\int_{T_k}g\ind{D}(0,\gx)
\Big(\sum_{n=(k-2)_+}^{k+2}\BBG_D[\mu_n](\gx)\Big)^q\,d\gx \label{upest3}
\end{align}
Then
\begin{equation}\label{upest0+}
  \BBG_D[(\BBG_D[\gm])^q](0)\leq 3^q(I_1+I_2+I_3)
\end{equation}
and we estimate each of the terms  on the right hand side separately.\2
\note{Estimate of $I_1$} We start with the following facts:
$$g\ind{D}(0,\gx)\leq c_N2^{k(N-2)} \forevery \gx\in T_k$$
and
$$g\ind{D}(\gx,z)\leq c_N2^{-n(2-N)} \forevery(\gx,z)\in T_k\ti F_n.$$
These inequalities and \req{mu3} imply, for every $\gx\in T_k$,
\[\BAL \BBG_D[\mu_n](\gx)&=\int_{F_n}g\ind{D}(\gx,z)d\mu_n(z)\leq c_N2^{n(N-2)}\mu_n(F_n)\\
&=c_N2^{n(N-2)}2^{-n(N-2q')}C_{2,q'}(2^nF_{n})=c_N2^{2n/(q-1)}C_{2,q'}(2^nF_{n}).\EAL\]
Hence
\begin{equation}\label{T3.1}
\BAL I_1&\leq c(N,q)\sum_{k=3}^\infty 2^{k(N-2)}\int_{T_k}\Big(\sum_{n=0}^{k-3}2^{2n/(q-1)}C_{2,q'}(2^nF_{n})\Big)^q\,d\gx\\
&\leq \sum_{k=3}^\infty 2^{k(N-2))}2^{-kN}\Big(\sum_{n=0}^{k-3}2^{2n/(q-1)}C_{2,q'}(2^nF_{n})\Big)^q\\
&\leq \sum_{k=3}^{M+1} 2^{-2k}\Big(\sum_{n=0}^{k-3}2^{2n/(q-1)}C_{2,q'}(2^nF^*_{n})\Big)^q.
\EAL\end{equation}
where $M=M(0)$ is defined as in \req{M(x)}.
Further, we claim that,
\begin{equation}\label{T3.2}\BAL
I'_1:=\sum_{k=3}^{M+1} 2^{-2k}\Big(\sum_{n=0}^{k-3}2^{2n/(q-1)}C_{2,q'}(2^nF^*_{n})\Big)^q&\leq\\
c(N,q)\sum_{n=0}^{M+1}2^{2n/(q-1)}C_{2,q'}(2^nF^*_{n}).
\EAL\end{equation}
 This  inequality is a \cons of the following statement proved in \cite[App. B]{MV3}:
\blemma{sum-int}
Let $K$ be a compact set in $\RN$ and let $\ga>0$ and $p>1$ be \sth $\ga p\leq N$. Put
\begin{equation}\label{gft}
  \gf(t)=C_{\ga,p}\big(\frac{1}{t}(K\cap B_{t})\big)= C_{\ga,p}\big(\rec{t}K\cap B_{1}\big) , \forevery t>0.
\end{equation}
Put $r_m=2^{-m}$. Then, for every  $\gg\in \BBR$ and every $k\in \BBN$,
\begin{equation}\label{sum-int}
\rec{c}\sum_{m=i+1}^{k}r_m^\gg\gf(r_m)
\leq \int_{r_k}^{r_i} t^\gg\gf(t)\frac{dt}{t}
\leq c\sum_{m=i+1}^{k}r_{m-1}^{\gg}\gf(r_{m-1}).
\end{equation}
where $c$ is a constant depending only on  $\gg,\,q,\,N$.
\es
\noindent Actually, in \cite{MV3} this result was proved in the case $\ga=2/q$, $p=q'$, in $\BBR^{N-1}$ assuming $2/(q-1)\leq N-1$. However
the proof
applies to any $\ga, p$ \sth $\ga p\leq N$. In particular it applies to the present case, namely, $\ga=2$, $ p=q'$ with $2q'\leq N$.

 We proceed to derive \req{T3.2} from the above lemma.
Put $r_m=2^{-m}$,  $\gg=-\frac{2}{q-1}$ and define $\gf$ and $\vgf$ by
\begin{equation}\label{a2.6i}
 \gf(r_{m}):=C_{2,q'}(r_m^{-1}F_m^*),\q \vgf(r,s):=\int_r^{s} t^{\gg}\gf(t)\frac{dt}{t}\q  0< r<s.
\end{equation}
By \rlemma{sum-int},
\begin{equation}\label{a2.6}
\rec{c}\sum_{m=i+1}^{k}r_{m}^{\gg}C_{2,q'}(r_m^{-1}F_m^*)\leq \vgf(r_k,r_i)\leq
 c\sum_{m=i+1}^{k}r_{m-1}^{\gg}C_{2,q'}(r_{m-1}^{-1}F_{m-1}^*),
\end{equation}
for every $i,k\in\BBN$, $i<k$. The constant $c$ depends only on $q,\, N,\, Q$.
Hence (taking into account that $F^*_m=\ems$ for $m>M+1$)
\begin{equation}\label{a2.6ii}\BAL
  \vgf(0,r_i)&:=\lim_{r \downarrow 0}\vgf(r,r_i)\leq c\sum_{m=i+1}^{\infty}r_{m-1}^{\gg}C_{2,q'}(r_{m-1}^{-1}F_{m-1}^*)\\
  &\leq c\sum_{m=i}^{M+1}r_{m}^{\gg}C_{2,q'}(r_{m}^{-1}F_{m}^*). 
\EAL\end{equation}
Further, by \req{a2.6},
\begin{align}\label{a2.7i}
I'_1=\sum_{k=3}^{M+1} r_{k}^2\Big(\sum_{n=0}^{k-3}r_n^\gg C_{2,q'}(r_n^{-1}F^*_{n})\Big)^q\leq
\sum_{k=3}^{M+1} r_{k}^2\vgf^q(r_{k-3},1).
\end{align}
Since $\vgf(\cdot,s)$ is non-increasing,
\begin{equation}\label{a2.7ii}
 \sum_{k=3}^{M+1} r_{k}^2\vgf^q(r_{k-3},1)\leq c\int_{r_{M-2}}^{1} t^{2}\vgf^q(t,1) \frac{dt}{t}\leq
 c\int_0^{1} t\vgf^q(t,1)\,dt.
\end{equation}
By  \req{a2.6ii}
\begin{align}\label{a2.7iii}
&\int_0^{1} t\vgf^q(t,1)\,dt\leq -c\int_0^{1}t^2\vgf^{q-1}(t,1)\dot{\vgf}(t,1) dt\\
&\leq -c\int_0^{1}\dot{\vgf}(t,1) dt\leq c \vgf(0,1)
\leq c\big(\sum_{m=0}^{M+1}r_{m}^{\gg}C_{2,q'}(r_{m}^{-1}F_{m}^*)\notag
\end{align}
Finally \req{a2.7i}--\req{a2.7iii} imply \req{T3.2}. In turn, \req{T3.1}, \req{T3.2} and \req{WF1} imply,
\begin{equation}\label{I1-fin}
 I_1\leq c(N,q)W_F(0).
\end{equation}

\noindent {\it Estimate of $I_2$}. Let $\gs>0$ and $\set{a_n}$ be a \seq of positive numbers.
Then,
$$\sum_{n=k}^\infty a_n\leq 2^{-\gs k}\Big(\frac{1}{1-2^{-\gs q'}}\Big)^{\rec{q'}}
\Big(\sum_{n=k}^\infty 2^{\gs nq}a_n^q\Big)^{\rec{q}}.$$
Applying this inequality with  $a_n=\BBG_D[\mu_n](\gx)$ we obtain
\begin{equation}\label{upest4'}\BAL
I_2&\leq c(N,q,\gs)\sum_{k=-1}^\infty \int_{T_k}g\ind{D}(0,\gx)2^{-\gs qk}
\sum_{n=k+2}^{\infty}2^{\gs nq}\BBG_D[\mu_n](\gx)^q\,d\gx\\
&\leq c\sum_{n\geq 1}2^{\gs nq}\sum_{1\leq k< n-2}{}\int_{T_{k}}2^{-\gs kq}g\ind{D}(0,\gx)
\BBG_D[\mu_n](\gx)^q d\xi\\
&\leq  c\sum_{n\geq 1}2^{\gs nq}\sum_{1\leq k< n-2}\int_{T_{k}}2^{-\gs kq}2^{k(N-2)}
\BBG_D[\mu_n](\gx)^q d\xi,
\EAL\end{equation}
where, in the last inequality, we used the fact that
$$g\ind{D}(0,\gx)\leq c_N2^{k(N-2)} \forevery \gx\in T_k.$$
Choosing $\gs=(N-1)/q$ we obtain,
\begin{equation}\label{upest4}
I_2\leq c(N,q)\sum_{n\geq 1}2^{n(N-1)}\sum_{1\leq k< n-2}\int_{T_{k}}2^{-k}\BBG_D[\mu_n](\gx)^q d\gx.
\end{equation}
Next we  estimate the term
\begin{equation}\label{Jkn}
  J_{k,n}:=\int_{T_{k}}\BBG_D[\mu_n](\gx)^q d\gx,
\end{equation}
in the case $1\leq k< n-2$. In view of \req{mu1} we have
$$ \BBG_D[\mu_n](\gx)=\int_{F_n}g\ind{D}(\gx,z)d\mu_n(z)=2^{-n(N-2q')}\int_{2^nF_n}\tl g(\gx', z')d\nu_n(z')$$
where $$\gx'=2^n\gx,\q \tl g(\gx', z')=g\ind{D}(2^{-n}\gx',2^{-n}z').$$
Observe that, if $\gx\in T_k$ then $\gx' \in T_{k-n}$. Thus
$$J_{k,n}=2^{-nN}\int_{T_{k-n}}\Big( 2^{-n(N-2q')}\int_{2^nF_n}\tl g(\gx', z')d\nu_n(z')
\Big)^q d\gx'.$$
Since
$$\tl g(\gx', z')\leq c_N2^{-n(2-N)}\absol{\gx'-z'}^{2-N}$$
we obtain
\begin{equation}\label{upest5}
 J_{k,n}\leq  c(N,q)2^{-n(N-2q')}\int_{T_{k-n}}\Big(\int_{2^nF_n}\absol{\gx'-z'}^{2-N}d\nu_n(z')
\Big)^q d\gx'.
\end{equation}
 Since $z'\in 2^nF_n\sbs B_1(0)$ while $\gx'\in T_{k-n}$, $k-n<-2$ it follows that $|\gx'|\geq 2$ and \consy
$$|\gx'-z'|\geq \rec{2}|\gx'|.$$
Therefore
\[\BAL &\int_{T_{k-n}}\Big(\int_{2^nF_n}\absol{\gx'-z'}^{2-N}d\nu_n(z')\Big)^q d\gx'
\leq c \nu_n(2^nF_n)^q\int_{T_{k-n}}|\gx'|^{(2-N)q}\,d\gx'\\
&\leq c(N,q)C_{2,q'}(2^nF_n)^q\int_{2^{n-k-1}}^{2^{n-k}}r^{(2-N)q+N-1}dr
\leq c(N,q)C_{2,q'}(2^nF_n) A(q,N)\EAL\]
where we used the fact that $C_{2,q'}(2^nF_n)\leq C_{2,q'}(B_1)$ and
$$A(q,N)=\begin{cases} 2^{(2-N)q+N}&\text{if }q>q_c\\ \ln 2 &\text{if }q=q_c.\end{cases}$$
Thus, for $k\geq n-2\geq -2$,
\begin{equation}\label{upest6}\BAL
 J_{k,n}&\leq  c(N,q)2^{-n(N-2q')}\norm{\nu_n}\ind{W^{-2,q}(\RN)}^q\\
 &=c(N,q)2^{-n(N-2q')}C_{2,q'}(2^nF_n).
\EAL\end{equation}
 By \req{upest4} and \req{upest6},
\begin{equation}\label{I2-fin}\BAL
 I_2&\leq c(N,q)\sum_{n\geq 0}2^{n(N-1)}\sum_{k\geq n-2}2^{-k}2^{-n(N-2q')}C_{2,q'}(2^nF_n)\\
&\leq c(N,q)\sum_{n\geq 0}2^{n(-2+2q')}C_{2,q'}(2^nF_n)\\
&=c(N,q)\sum_{n\geq 0}2^{\frac{2n}{q-1}}C_{2,q'}(2^nF_n)=c(N,q)W_F(0).
\EAL\end{equation}

\noindent {\it Estimate of $I_3$}\hskip 2mm By \req{upest3} and the notation \req{Jkn} we have
\begin{equation}\label{I3-1}\BAL
 I_3&\leq 5^q\sum_{k=-1}^\infty\int_{T_k}g\ind{D}(0,\gx)
\sum_{n=(k-2)_+}^{k+2}\BBG_D[\mu_n](\gx)^q\,d\gx\\
&\leq \sum_{k=-1}^\infty 2^{k(N-2)}\sum_{n=(k-2)_+}^{k+2}J_{k,n}.
\EAL\end{equation}
By \req{upest5}
$$J_{k,n}\leq  c(N,q)2^{-n(N-2q')}\int_{T_{k-n}}\Big(\int_{2^nF_n}\absol{\gx'-z'}^{2-N}d\nu_n(z')\Big)$$
and, in the present case $-2\leq n-k\leq 2$. Therefore $T_{k-n}\sbs B_{4}(0)$ and \consy, for $(\gx',z')$ in the domain of integration
of the integral above,
$$\absol{\gx'-z'}^{2-N}\approx \CB_2(\gx',z')$$
where $\CB_2$ denotes the Bessel kernel with index 2. Hence,
\[\BAL  &\int_{T_{k-n}}\Big(\int_{2^nF_n}\absol{\gx'-z'}^{2-N}d\nu_n(z')
\Big)^q d\gx'\leq\\
&c(N,q)\norm{\nu_n}\ind{W^{-2,q}(\RN)}^q=c(N,q)C_{2,q'}(2^nF_n).
\EAL\]
Therefore,
\begin{equation}\label{I3-fin}\BAL
I_3&\leq c(N,q)\sum_{k=-1}^\infty 2^{k(N-2)}\sum_{n=(k-2)_+}^{k+2}2^{-n(N-2q')}C_{2,q'}(2^nF_n)\\
&\leq c(N,q)\sum_{k=-1}^\infty 2^{k(N-2)}2^{-k(N-2q')}C_{2,q'}(2^kF_k)\\
&=c(N,q)\sum_{k=-1}^\infty 2^{2k/(q-1)}C_{2,q'}(2^kF_k)\leq c(N,q)W_F(0)
\EAL\end{equation}

Combining \req{upest0+} with the inequalities \req{I1-fin}, \req{I2-fin} and \req{I3-fin} we obtain
\begin{equation}\label{Gmuq}
      \BBG_D[(\BBG_D[\gm])^q](0)\leq c(N,q)W_F(0).
\end{equation}
Finally, we combine \req{Gmu} with \req{G[mu]} and \req{Gmuq} and  replace $\mu$ by $\ge\mu$, $\ge>0$,  to obtain
\begin{equation}\label{cap-fin}
  u^{\ge\mu}(0)\geq (c_1(N,q)\ge -c_2(N,q)\ge^q)W_F(0). 
\end{equation}
Choosing $\ge:=(c_1(N,q)/2c_2(N,q))^{1/(q-1}$
we obtain \req{WF<u;D} with $c(N,q)=c_1(N,q)\ge/2$.
\eproof

\section{Properties of $U_F$ for $F$ compact}
As before we assume that $F$ is a compact set. Combining the capacitary
estimates contained in Theorems \ref{t:up-est}, \ref{t:low-est} and \req{WF1}we have
\begin{equation}\label{UsimW}
 U_F\sim W_F \sim W^*_F \q \text{in } D=\RN\sms F
\end{equation}
In the present section we use this result in order to establish
several properties of the maximal solution.
\subsection{\normalfont{\it The maximal solution is $\gs$-moderate}}
\bth{gsmod} $U_F=V_F\,$; \consy $U_F$ is \gsmod.
\es
\bproof By \req{UsimW} there exists a constant  $c=c(N,q)$ \sth
\begin{equation}\label{cons1.2}
 U_F\leq cV_F\,.
\end{equation}
 If the two solutions are not identical we have
\begin{equation}\label{u<U}
V_F\,(x)<U_F(x) \forevery x\in D.
\ee
Let $\ga=\rec{2c}$ and put $v=(1+\ga)V_F\,(x)-\ga U_F$. Then $\ga V_F\,(x)<v<U_F$ and (as $0<\ga<1$) $\ga V_F\,(x)$
is a subsolution of \req{eq} in $D$. As in
\cite{MV1} we find that $v$ is a supersolution. It
follows that there exists a solution $w$ \sth $\ga V_F\,(x)\leq w\leq v<V_F\,(x)$. But, by the definition of
$V_F\,$ (see \req{ca4}), it is easy to see that the
smallest solution of \req{eq} dominating $\ga V_F\,(x)$ is $V_F\,(x)$. Therefore $w=V_F\,(x)$. This contradicts
\req{u<U}.

By a standard argument, the definition of $V_F\,(x)$ implies that it is \gsmod.
\eproof
\subsection{\normalfont{\it A continuity property of $U_F$ relative to capacity}}
\blemma{CKsmall}
There exists a positive constant $c$ depending only on $N,q$ \sth, for every compact set $K\sbs B_1(0)$,
there exists an open \ngh $N_K$ of $K$ \sth
\begin{equation}\label{CK1}
 C_{2,q'}(N_K)\leq 4C_{2,q'}(K) \q\text{\rm and }\; \int_{B_1(0)\sms N_K}U_K\,dx\leq c C_{2,q'}(K).
\ee
\es
\note{Note} In general $\int_{B_1(0)\sms K}U_K\,dx$ may be infinite.
Of course, \req{CK1} is meaningful only if  $4C_{2,q'}(K)<C_{2,q'}(B_1(0))$.

\bproof  Let  $ \bar c$ be the constant in \req{UK1} with $R=2$. Assume that
\begin{equation}\label{CK2}
\C2q(K)\leq a:=\C2q(B_1)/8
\ee
 and pick $\gg_1$ so that
\begin{equation}\label{CK4}
0<\gg_1\leq \C2q(K).
\ee
By \rlemma{intUK} and \req{cap-def} there exists $\eta\in X_K(\RN)$ \sth
\begin{equation}\label{CK5}
\BAL\norm{\eta}^{q'}\ind{W^{2,q'}(\RN)}&\leq C_{2,q'}(K)+\gg_1,\\
 \int_{B_2(0)\sms K}U_K(1-\eta)^{2q'}\,dx&\leq \bar c(C_{2,q'}(K)+\gg_1).
 \EAL\ee
Fix $\eta$ and denote,
$$K_{(\ga)}=\set{x\in B_1(0): (1-\ga)\leq \eta} \forevery \ga\in (0,1).$$
Then $K\sbs K_{(\ga)}$ and
\[\BAL C_{2,q'}(K_{(\ga)})&\leq (1-\ga)^{-q'}\norm{\eta}^{q'}\ind{W^{2,q'}(\RN)}\\
&\leq (1-\ga)^{-q'}(C_{2,q'}(K)+\gg_1) \leq \frac{2\C2q(K)}{(1-\ga)^{q'}}.
\EAL\]
Therefore, using \req{CK2}, we obtain
\[ (1-\ga)^{-q'}=2 \Lra C_{2,q'}(K_{(\ga)})\leq 4\C2q(K)\leq \C2q(B_1)/2.\]
Hence, by \req{CK5},
\begin{equation}\label{CK6}
\int_{B_2(0)\sms K_{(\ga)}}U_K\,dx\leq \bar c \ga^{-2q'}(C_{2,q'}(K)+\gg_1)\leq (4\bar c)C_{2,q'}(K)
\end{equation}
where $\ga=1-2^{-1/q'}$.
\eproof

\subsection{\normalfont{\it Wiener criterion for blow up of $U_F$}}

\bth{Wiener} For every point $y\in F$,
\begin{equation}\label{Wiener}
 \lim_{F^c\ni x\to y} U_F(x)= \infty \iff W_F(y)=\infty.
\end{equation}
\es
\bproof
Without loss of generality we may assume that $y=0$ and that $F\sbs B_1(0)$. In order to justify the second part of this remark
we observe that, for every $m\in\BBN$,
\begin{equation}\label{similarity}\BAL
  2^{-2m/(q-1)}U_{F}(2^{-m} x)&=U\ind{2^mF}(x) \forevery x\in (2^mF)^c,\\  W_{F}(0)&=2^{2m/(q-1)}W\ind{2^mF}(0).
\EAL\end{equation}
Denote
\begin{equation}\label{am(y)}
 a_m(x)=C_{2,q'}\left(2^{m}F_m(x)\right), \q a^*_m(x)=C_{2,q'}\left(2^{m}F^*_m(x)\right)
\end{equation}
There exists a constant $c=c(N,q)$ \sth for every Borel set $A\sbs B_1(0)$,
\begin{equation}\label{Wiener1}
 C_{2,q'}(2A)\leq cC_{2,q'}(A).
\end{equation}
If $x,\gx\in \RN$, $|x-\gx|\leq r_m=2^{-m}$ and $0\leq k\leq m$ then
$$2^{k}F^*_k(\gx)=2^k(F\cap B_{r_k}(\gx))\sbs 2^k(F\cap B_{2r_k}(x))=2(2^{k-1}(F\cap B_{r_{k-1}}(x)).$$
Hence
\begin{equation}\label{Wiener2}
 a^*_k(\gx)\leq ca^*_{k-1}(x) \text{ for }0\leq k\leq m, \q \sum_{k=0}^m 2^{2k/(q-1)}a^*_k(\gx)\leq cW^*_F(x).
\end{equation}
As $x\to \gx$, $m\tin$ and we obtain
\begin{equation}\label{W<climinf}
W^*_F(\gx)\leq c(N,q)\liminf_{x\to\gx}W^*_F(x).
\end{equation}
By \req{UsimW}, \req{W<climinf} implies that
\req{Wiener} holds in the direction $\Lla\,$.

In order to prove \req{Wiener} in the opposite direction
we derive the inequality,
\begin{equation}\label{cW>liminf}
  \liminf_{F^c\ni x\to0}W^*_F(x)\leq c(N,q)W^*_F(0).
\ee
If $W^*_F(0)=\infty$ there is nothing to prove. Therefore we assume that

$$W^*_F(0)= \sum_{-\infty}^{M(0)} 2^{\frac{2m}{q-1}}C_{2,q'}\left(2^{m}F^*_m(0)\right)<\infty.$$
By \rlemma{CKsmall}, with $K_m=2^{m}F^*_m(0)$, there exists an open \ngh $G_m$ of $K_m$ \sth
$$ \C2q(G_m)\leq 4 \C2q(K_m) \txt{\rm and} \int_{B_1(0)\sms G_m}U_{K_m}dx\leq c(N,q)\C2q(K_m).$$

 Put $T':=[5/8\leq |x|\leq 7/8]$ and let
$E_m$ be a compact subset of $T'\sms G_m$  \sth
$$\BAL
\C2q(E_m)&>\rec{2}\C2q(T'\sms G_m)
\geq \rec{2}\C2q(T')-2\C2q(K_m)\\
&\geq \rec{2}\C2q(T')-2^{1-\frac{2m}{q-1}}W^*_F(0).\EAL$$
Therefore, there exists an integer $m_0$ \sth, for $m\geq m_0$,
\begin{equation}\label{cW1}\BAL
  \inf_{E_m} U_{K_m}&\leq |E_m|^{-1}\int_{E_m}U_{K_m} dx\\
  &\leq |E_m|^{-1}c(N,q)\C2q(K_m)
  \leq c(N,q)\C2q(K_m)/\C2q(E_m)\\
  &<4c(N,q)\C2q(B_1(0))^{-1}\C2q(K_m)=c(N,q)\C2q(K_m).
\EAL\end{equation}
Hence, by \req{similarity},
$$\inf_{2^{-m}E_m}U_{F^*_m(0)}=2^{2m/(q-1)} \inf_{E_m} U_{K_m}<c(N,q)2^{2m/(q-1)}\C2q(K_m),$$
which implies, for $m\geq m_0$,
\begin{equation}\label{cW2}
\inf_{(2^{-m}T')\sms F}U_{F^*_m(0)}\leq c(N,q)2^{2m/(q-1)}\C2q(K_m)\leq c(N,q)W^*_F(0).
\ee

\par Fix $j> m_0$ and let $\gx\in (2^{-j}T')\sms F$. Denote
$$F^j:=F\sms B_{2^{-j+1}}(0),\q E_{k,j}:=F^j\cap \set{x\in B_1(0):|x-\gx|\leq 2^{-k}},\q \bar j:=\left[\frac{j}{8}\right].$$
Since $\dist(\gx,F^j)\geq 2^{-j/8}$,
$$W_{F^j}(\gx)=\sum_{-\infty}^{\bar j}2^{2k/(q-1)}\C2q(2^k E_{k,j}).$$
For $k\leq \bar j$
$$x\in E_{k,j}\Lra |x|\leq 2^{-k}+2^{-j}\leq 2^{-k}+2^{-8(k-1)}\leq 2^9 2^{-k}.$$
Thus $E_{k,j}\sbs F^*_{k-9}$ and
\begin{equation}\label{cW3}
W_{F^j}(\gx)=\sum_{-\infty}^{\bar j}2^{2k/(q-1)}\C2q(2^kF^*_{k-9}(0))\leq c(N,q)W^*_F(0)
\end{equation}
for every $\gx\in (2^{-j}T')\sms F$.
By \req{cW2}, we can choose $\gx^{j}\in (2^{-j}T')\sms F$ \sth
$$U_{F^*_j(0)}(\gx^j)\leq c(N,q)W^*_F(0).$$
Hence, by \req{cW3}, bearing in mind that $U_K\sim W_K\sim W^*_K$ for every compact $K$ we obtain
\begin{equation}\label{cW4}
U_F(\gx^j)\leq U_{F^*_j(0)}(\gx^j)+ U_{F^j}(\gx^j)\leq c(N,q)W^*_F(0) \forevery j\geq m_0.
\ee
This implies \req{cW>liminf} and completes the proof.
\eproof
\bcor{tlWF} Define
$$\wtl W_F(x)=\liminf_{y\to x}W_F(y) \forevery x\in \RN.$$
Then,  $\tl W_F$ is l.s.c. in $\RN$ and $\wtl W_F\sim W_F$.
In addition,  $\tl W_F$  satisfies Harnack's inequality in compact subsets of $\RN\sms F$.
\es
\bproof The lower semi-continuity of $\wtl W_F$ follows from its definition. The equivalence $\wtl W_F\sim W_F$
 follows from  \req{W<climinf} and \req{cW>liminf}. The last statement follows from the fact that
$U_F$ satisfies Harnack's inequality and $\tl W_F\sim U_F$.
\eproof
\vskip 2mm
\subsection{\normalfont{\it $U_F$ is an almost large solution}}
\bth{almostlarge} For every compact set $F\sbs \RN$, $U_F$ is an almost large solution.
\es
\bproof In view of \rth{Wiener} it is enough to show that there exists a set $A\sbs F$ \sth
\begin{equation}\label{almost1}
 \C2q(A)=0\txt{and}  W_F(y)=\infty \forevery y\in F\sms A.
\end{equation}
It is known (see \cite[Ch. 6]{AH}) that every point in $F$, with the possible exception of a set $A_1$ of $\C2q$ capacity zero,
is a $\C2q$-thick point of $F$, i.e.,
\begin{equation}\label{thkpoint}
  \Lambda^{2,q'}_F(y):=\sum_{0}^{\infty}\Big( 2^{2m/(q-1)}\C2q(F^*_m(y))\Big)^{q-1}=\infty \forevery y\in F\sms A_1.
\end{equation}
We show that,
\begin{equation}\label{thk-WF}
    \Lambda^{2,q'}_F(y)\leq c(N,q)(W_F(y))^{\tl q}, \q \tl q=\min(1,q-1) \forevery q>1.
\end{equation}

Recall that,
$$\C2q(F^*_m(y))\leq c(N,q)2^{-2m/(q-1)}\C2q(2^mF^*_m(y))$$
so that
\begin{equation}\label{thk1}
\Lambda^{2,q'}_F(y)\leq c(N,q)\sum_{0}^{\infty}\big(\C2q(2^mF^*_m(y))\big)^{q-1}.
\ee
 In view of the fact that $\C2q(2^mF^*_m(y))\leq \C2q(B_1)$, if $q\geq 2$, \req{thk1} implies \req{thk-WF}.
If $1<q<2$,
$$\BAL &\sum_{0}^{\infty}\big(\C2q(2^mF^*_m(y))\big)^{q-1}\leq\\
\Big( &\sum_{0}^{\infty} 2^{-\frac{2m(2-q)}{q-1}} \Big)^{2-q}\Big( \sum_{0}^{\infty}2^{\mq} \C2q(2^mF^*_m(y)) \Big)^{q-1},
\EAL $$
which again implies \req{thk-WF}. Clearly \req{thkpoint} and \req{thk-WF} imply \req{almost1}.
\eproof
\section{'Maximal solutions' on arbitrary sets and uniqueness I}\label{arbitrary}

For any Borel set $E$ put
\begin{align}
  \CT(E)&:=\set{\mu\in W^{-2,q}_+(\RN): \;\mu(E^c)=0},\label{CTE}\\
  V_E&:=\sup\set{u_\mu:\mu\in \CT(E)},\label{VE}
\end{align}
where $u_\mu$ denotes the solution of \req{equ-meas}. If $\C2q(E)=0$ the only  measure $\mu\in W^{-2,q}_+(\RN)$ that is concentrated
on $E$ is the measure zero. Therefore in this case $V_E=0$.

By \rth{gsmod},
\begin{equation}\label{Ecompact}
 E\text{ compact}\;\Lra\; V_E=U_E.
\end{equation}
 Therefore the definition of $V_E$ can be seen as an extension,  to general sets,
of the notion of 'maximal solution', previously defined for compact sets. However, by its definition, $V_E$
{\em dominates only
\gsmod solutions in $E^c$,} i.e., solutions of the form $\lim u_{\mu_n}$ where $\set{\mu_n}$ is an increasing \seq
of measures in $W^{-2,q}_+(\RN)$ concentrated in $E$.

At this stage, it is not clear in which sense $V_E$ is a solution of \req{eq}
in $E^c$, which, in general, is not an open set.
 This question will be
discussed in the following sections.

The $\C2q$ fine topology (see \cite{AH} for definition and details) plays a central role in
the remaining part of the paper.
 If $A$ is a set in $\RN$ we
denote by $\wtl A$ the closure of $A$ in the $\C2q$ fine topology and
by $\mathrm{int}_q A$  the interior of $A$ relative to this
topology.

Recall that a set $A\sbs \RN$ is {\em $\C2q$-quasi open} if, for every $\ge>0$, there exists an open set $G_\ge$ \sth
$$A\sbs G_\ge,\qq \C2q(G_\ge\sms A)<\ge.$$
A set is {\em $\C2q$-quasi closed} if its complement is quasi open.

Every \qfine open set is $\C2q$-quasi-open. On the other hand, if $E$ is $\C2q$-quasi open then
(see \cite[Section 6.4]{AH})
$$\C2q(E\sms \qinterior E)=0.$$

This implies that every $\C2q$-quasi closed set $F$ can be written in the form
$$F=\bigcup_n^\infty K_n\bigcup Z,$$
where $\set{K_n}$ is an increasing \seq of compact sets and
$$\C2q(F\sms K_n)\to 0,\qq \C2q(Z)=0.$$
Furthermore, if $E$ is $\C2q$-quasi closed then
$$\C2q(\wtl E\sms E)=0.$$

In the first two theorems below we describe some basic properties of $V_E$. These results are then used in order
to establish
a rather general uniqueness result for almost large solutions.

\bth{VFest} Let $F$ be a  \quasi closed set. Then
\begin{equation}\label{VFtin}
  \lim_{F^c \ni x\to y}V_F(x)=\infty \q \text{for $\C2q$-a.e. } y\in F
\end{equation}
and $V_F$ satisfies
\begin{equation}\label{W<V<cW}
  \rec{c}W_F\leq V_F\leq cW_F,
\end{equation}
where $c$ depends only on $N,q$.                   
Finally, for every $x\in F^c$,
\begin{equation}\label{VEto0}
 W_F(x)<\infty\; \Lra \;\lim_{\BSA{c} \C2q(E)\to0\\ E\sbs F\ESA}V_E(x)=0.
 \ee
\es
\bproof There exists an increasing \seq of  compact sets $\set{K_n}$ \sth $K_n\sbs F$ and $\C2q(F\sms K_n)\to0$.
 By \rth{gsmod} $U_{K_n}=V_{K_n}$ and, obviously, $ V_{K_n}\leq V_F$.  By
 \rth{almostlarge}, \req{VFtin} holds if $F$ is replaced by $K_n$. Therefore, by taking the limit as $n\tin$, we obtain
  \req{VFtin} in the general case.

 If $\mu\in \CT_F$ then $u_\mu=\lim u_{\mu_n}$ where $\mu_n=\mu\chr{K_n}$. Therefore
\begin{equation}\label{VF=lim}
 V_F=\lim U_{K_n}.
\end{equation}
Since $U_{K_n}$ satisfies estimates  \req{UF<WF} and \req{WF<u} for every $n$, it follows that $V_F$ satisfies \req{W<V<cW}.

We turn to the proof of the last assertion. Let $\set{E_j}$ be a \seq of subsets of $F$ \sth $\C2q(E_j)\to 0$.
We must show that
\begin{equation}\label{WF=infty}
 \gx\in F^c,\q  \limsup V_{E_j}(\gx)>0 \;\Lra \; W_F(\gx)=\infty.
\end{equation}

By taking a \sseq we may assume that there exists $a>0$ \sth $V_{E_j}(\gx)>a$ for all $j$.
Since $V_{\wtl E_j}(\gx)=V_{E_j}(\gx)$ and $\C2q(E_j)\to 0$ implies $\C2q(\wtl E_j)\to 0$ we may assume that the sets $E_j$ are
\qfine closed.
By \req{VF=lim}  it follows
that, for every $j$, there exists a compact set $K_j\sbs E_j$ \sth
\begin{equation}\label{UKj>a}
  U_{K_j}(\gx)>a.
\end{equation}

\bcom
Given $R>0$ put
$K_{j,R}=K_j\cap \bar B_R(\gx)$ and $K'_{j,R}=K_j\sms B_R(\gx)$. Then
$$V_{K'_{j,R}}(\gx)\leq c(N,q)R_{2/(j-1)},$$
so that choosing $R $ sufficiently large we obtain,
$$V_{K_{j,R}}(\gx)\geq V_{K_{j}}(\gx)- V_{K'_{j,R}}(\gx)\geq a/2.$$
 Thus there exists a \seq of {\em compact sets} $\set{K_j}$ \sth
 $$S_j\sbs B_R(\gx)\cap F,\q  \C2q(S_j)\to 0,\q V_{S_j}(\gx)>a/2 \txt{for $j=1,2\ldots$.}$$
\end{comment}

By negation, suppose that $W_F(\gx)<\infty$. Then
$$\lim_{J\tin}\sum_J^\infty 2^{2j/(q-1)}\C2q(2^jF_j(\gx))\to 0,$$
$F_j$ being defined as in \req{Fm}. Pick a positive integer $J$ \sth
\begin{equation}\label{partofWF}
\sum_J^\infty 2^{2j/(q-1)}\C2q(2^jF_j(\gx))<a/4C,
\end{equation}
where $C$ is the constant in \req{UF<WF}.

Pick a \sseq  of $\set{K_j}$, say $\set{K_{j_n}}$, \sth $\C2q(K_{j_n})<\ge/2^n$, with $\ge$ to be determined.
The set   $A:=\bigcup_1^\infty K_{j_n}$ is \quasi closed and $\C2q(A)\leq \sum_1^\infty\C2q(K_{j_n})<\ge$. Further,
$$\BAL
W_A(\gx)&=\sum_{-\infty}^\infty 2^{2j/(q-1)}\C2q(2^jA_j(\gx))\leq  \sum_{-\infty}^{-1}2^{2j/(q-1)}\C2q(2^jA)\\
&+\sum_{0}^{J-1} 2^{2j/(q-1)}\C2q(2^jA)+\sum_{J}^\infty 2^{2j/(q-1)}\C2q(2^jF_j(\gx)),
\EAL$$
where $A_j(\gx)$ is defined as in \req{Fm} with $F$ replaced by $A$. (We  used the fact that $A_j(\gx)\sbs A\sbs F$.)
By \req{C(aE),1},
$$\sum_{-\infty}^{-1}2^{2j/(q-1)}\C2q(2^jA)\leq c_1(N,q) \sum_{-\infty}^{-1} 2^{jN}\C2q(A)\leq c_2(N,q)\ge.$$
By \req{C(aE),2},
$$\sum_{0}^{J-1} 2^{2j/(q-1)}\C2q(2^jA)\leq c_1(N,q, J)\sum_{0}^{J-1} 2^{2jN}\C2q(A)\leq c_2(N,q,J)\ge.$$
Therefore,   choosing $\ge=(a/4C)(c_2(N,q)+c_2(N,q,J))^{-1}$ and using \req{partofWF} we obtain,
$$W_A(\gx)<a/2C.$$
Since $V_A$ satisfies \req{UF<WF} we conclude that $V_A(\gx)<a/2$.
As $$U_{K_{j_n}}=V_{K_{j_n}}\leq V_A,$$ this contradicts \req{UKj>a}.
\eproof

\bth{V=limU} Let $F$ be a \quasi closed set and let $\set{F_n}$ be an increasing \seq of compact subsets of $F$
\sth $\C2q(F\sms F_n)\to 0$. Then
\begin{equation}\label{V=limU}
  V_F=\lim U_{F_n}.
\end{equation}
Furthermore, there exists an increasing \seq of non-negative measures $\set{\mu_n}\sbs W^{-2,q}(\RN)$ \sth
$\mu_n(F_n^c)=0$ and $u_{\mu_n}\to V_F$ in $F^c$.

Finally, for every $y\in F$,
\begin{equation}\label{WF-qclosed}
 W_F(y)<\infty\; \Llra \;\liminf_{F^c\ni x\to y}V_F(x)<\infty.
 \ee
\es
\bproof From the definition of $V_F$ it follows that $V_F=\lim V_{F_n}$. By \rth{gsmod}, $V_{F_n}=U_{F_n}$.

Let $\gx\in D_n=F_n^c$ and let $\set{\tau^n_k}_{k=1}^\infty$ be a \seq in $W^{-2,q}(\RN)$ \sth $\tau^n_k(D_n)=0$
and $u\ind{\tau^n_k}(\gx)\to U_{F_n}(\gx)$. Note that $w^n_m=\max(u\ind{\tau^n_1},\ldots, u\ind{\tau^n_m})$ is a subsolution
of the equation
$$-\Gd w+ w^q=\mu^n_m:=\max({\tau^n_1},\ldots, {\tau^n_m}) \txin{D_n}.$$
Therefore $v^n_m=u\ind{\mu^n_m}$ is the smallest solution in $D_n$ dominating $w^n_m$.
The \seq $\set{v^n_m}_{m=1}^\infty$ is increasing, bounded by $U_{F_n}$
and
$v^n:=\lim_{m\tin}v^n_m$ is a solution of \req{eq} in $D_n$ \sth $v^n(\gx)=U_{F_n}(\gx)$. The fact that $v^n\leq U_{F_n}$ and equals it
at a point $\gx\in D_n$ implies that $v^n= U_{F_n}$.

Put
\begin{equation}\label{taun}
\tau^{(n)}:=\sum_m a^n_m\mu^n_m,\q a^n_m:=2^{-m}\norm{\mu^n_m}_{W^{-2,q}_+(\RN)}.
\end{equation}
Then
\begin{equation}\label{Fn-taun}
U_{F_n}=\lim_{k\tin}u\ind{k\tau^{(n)}}.
\ee
Finally, if $\mu_n:=\sum_1^n\tau^{(j)}$ then $\set{\mu_n}$ is increasing and $u_{\mu_n}\to V_F$.

The last statement of the theorem is proved exactly as in the case that $F$ is compact (see \rth{Wiener}).
\bcom
\begin{equation}\label{tau}
\tau:=\sum_n b_n\tau^n,\q b_n:=2^{-n}\norm{\tau^n}_{W^{-2,q}_+(\RN)}
\ee
and observe that
\begin{equation}\label{F-tau}
V_F=\lim u_{k\tau}.
\ee
\end{comment}
\eproof

\bth{VtlE} Let $E$ be a Borel set \sth $\C2q(E)>0$. Then
\begin{equation}\label{VE1}
   V_E=V_{\wtl E}
\end{equation}
and, if $ \mu\in W^{-2,q}_+(\RN)$,
\begin{equation}\label{VE2}
u_\mu<V_E \iff \mu(\RN\sms \tl E)=0.
\end{equation}
Furthermore, if $E$ is \quasi closed, there exists
$\tau\in W^{-2,q}_+(\RN)$ \sth $\tau(E^c)=0$ and
\begin{equation}\label{VE3}
V_E=\lim_{k\tin}u_{k\tau}.
\ee
\es
\bproof
We prove \req{VE1} under the assumption that $E$ is bounded, say,  $\ovl E\sbs B_R$.
For the general case we observe  that
$$\lim_{R\tin}V_{E\cap B_R}= V_E.$$
\underline{Assertion 1}: {\em Let $\ovl\CT(E)$ denote the closure of $\CT(E)$ in $W^{-2,q}_+(B_R)$. If $\mu\in\ovl\CT(E)$
then $u_\mu\leq V_E$.}\2
Let $\set{\nu_n}$ be a \seq in $\CT(E)$ \sth $\nu_n\to \nu$ in $W^{-2,q}_+(B_R)$. Let $u_n$ be the solution of
$$-\Gd u_n+u_n^q=\nu_n \txin{B_R},\q u_n=0\txon{\Gs_R}.$$
Then $\set{u_n}$ converges in $L^q(B_R)$ and the limit $u$ is a weak solution of
$$-\Gd u+u^q=\nu \txin{B_R},\q u=0\txon{\Gs_R}.$$
Since  $u_n\leq V_E$  it follows that $u=u_\nu \leq V_E$. \2
\underline{Assertion 2}:
\begin{equation}\label{barCTE}
 \nu\in\CT( \wtl E)\Lra \nu\in \ovl{\CT}(E).
\end{equation}
 Suppose that $\nu\in\CT( \wtl E)$ but $\nu\not\in \ovl{\CT}_E$.
Then there exists $\gf\in W^{2,q'}(\RN)$ \sth
\begin{equation}\label{neg1}
\norm{\gf}_{W^{2,q'}(\RN)}=1,\q \langle \gf, \nu\rangle>0, \q \langle \gf, \mu\rangle=0
\forevery \mu\in{\CT}_E.
\end{equation}
We choose $\gf$ to be a \qfine continuous representative of its equivalence class
(see \cite[Proposition 6.1.2]{AH}).
Thus the inverse image (by $\gf$) of every open interval is quasi-open (see \cite[Proposition6.4.10]{AH}).
It follows that $$A_0:=\set{\gs:\,\gf(\gs)=0} \txt{is \qfine closed.}$$
 We  show that
\begin{equation}\label{neg2}
\C2q(\wtl E\sms A_0)=0.                         
\end{equation}
Put $A_1:=\wtl E\sms A_0$ and
$$A_1^+=\set{x\in A_1:\,\gf(x)>0},\q A_1^-=\set{x\in A_1:\,\gf(x)<0}.$$
If \req{neg2} does {\em not} hold then
$$\text{either}\q \C2q(A_1^+)>0, \txt{or} \C2q(A_1^-)>0.$$
Each of these sets is \qfine open relative to $\wtl E$, i.e., there exist \qfine open sets $Q_1,Q_2$ \sth
$Q_1\cap \wtl E=A_1^+$ and
$Q_2\cap \wtl E=A_1^-$. If, say, $\C2q(Q_1\cap \wtl E)>0$ then $\C2q(Q_1\cap E)>0$
(because $\C2q(G)\sim \C2q(\wtl G)$ for any Borel set $G\sbs B_R$).
Let $\mu\in W^{-2,q}_+(\RN)$ be a non-trivial measure ,  supported in a compact subset of $Q_1\cap E$. Then
$$\langle \gf, \mu\rangle>0.$$
This contradicts \req{neg1} and proves \req{neg2}.

Further \req{neg2} implies that $\gf=0$ $\C2q$-a.e. on $\wtl E$ which implies $\langle \gf, \nu\rangle=0$
in contradiction to \req{neg1}. This proves Assertion 2.\1
\bcom
{\em Proof of \req{barCTE}  $\Lla$.} Suppose that $\nu\in \ovl{\CT}(E)$ but $\nu\not \in \CT(\tl E)$.
Then there exists a compact set $A\sbs \prt_q E$ \sth $\nu(A)>0$. On the other hand,
$$\langle \gf,\nu\rangle=0 \forevery \gf \in W^{2,q'}(\RN) \;\text{ \sth $\gf=0$ on $\tl E$}.$$
To verify this,
let $\set{\nu_n}$ be
a \seq in $\CT(E)$ \sth $\nu_n\to\nu$ in $W^{-2,q}(\RN)$ and let
 $\gf\in W^{2,q'}(\RN)$ be a function \sth $\gf=0$ on $E$. Then $\langle\gf,\nu_n\rangle=0$ and \consy
$\langle\gf,\nu\rangle=0$. But $\gf=0$ on $E$ \ifif $\gf=0$ on $\tl E$. Therefore $\nu(\RN\sms \tl E)=0$. \2
\underline{Assertion 3}:  $\ovl{\CT}(E)=\CT(\tl E)$.\2
See Lemma 1.12 in NOTE23.\\
\end{comment}
\indent Combining these assertions we conclude:
\begin{equation}\label{pos1}
\nu\in \CT(\wtl E)\Lra u_\nu\leq V_E\Lra V_{\wtl E}=sup\set{u_\nu:\,\nu\in \CT(\tl E)}\leq V_E.
\end{equation}
Since, trivially, $V_E\leq V_{\wtl E}$ we obtain \req{VE1}.

If $E$ is \quasi closed then, by \rth{V=limU}, there exists an increasing \seq $\set{\mu_n}$ in
$W^{-2,q}_+(\RN)$ \sth $\mu_n(E^c)=0$ and $u_{\mu_n}\to V_E$. Put
\begin{equation}\label{VE3'}
\tau:=\sum a_n\mu_n,\q a_n:=2^{-n}\norm{\mu_n}_{W^{-2,q}_+(\RN)}.
\end{equation}
Then $\tau\in W^{-2,q}_+(\RN)$, $\tau(E^c)=0$ and \req{VE3} holds.

We turn to the proof of \req{VE2}. The implication $$  \mu(\RN\sms \tl E)=0\Lra u_\mu<V_E$$
is a \cons of \req{VE1}. To prove the implication in the opposite direction we may assume
that $E$ is compact. (This follows from \rth{V=limU}.) By negation,
 suppose there exists $\mu\in W^{-2,q}(\RN)$ such that $u_\mu<V_E$ but $\mu(\RN\sms \tl E)>0$.
 It follows that there exists a compact set $K\sbs \RN\sms \tl E$ \sth
 $\mu(K)>0$. Let $v_n:=u\ind{n\mu\chr{K}}$. Then $v_n\leq n u_\mu$ because $n u_\mu$ is a supersolution of
 the equation $-\Gd w+w^q=n\mu\chr{K}$. On the other hand, $V_E$ is the largest solution dominated by $nV_E$, for every $n$.
 Therefore
 \begin{equation}\label{v<VE}
  v=\lim v_n\leq V_E.
 \end{equation}

If $A$ is an  open \ngh of $K$ \sth $\dist(A,E)>0$ then
 $V_E\in L^q(A)$. On the other hand
 $$\int_{A\sms K}v^q=\infty.$$
 Therefore $(v-V_E)_+$ is positive in an open subset of $A\sms K$.
 This contradicts \req{v<VE}.
\eproof

\bcom \bth{uniq} Let $\Gw$ be a domain with compact boundary in
$\RN$. Put $F=\prt\Gw$ and
 $D:=\RN\sms \ovl \Gw$.
Assume that $\C2q$-a.e. point $x\in F$ is a $\C2q$-thick point of $D$, i.e.,
\begin{equation}\label{uniq1}
    \C2q(F\sms \wtl D)=0.
\end{equation}
 If $v$ is a non-negative solution of the \bvp
\begin{equation}\label{bvp-infty}
  -\Gd v+ v^q=0 \txin{ \Gw},\q \lim_{\Gw \ni x\to y}v(x)=\infty \q \text{$\C2q$-a.e. } y\in F
\end{equation}
then $v=U_F$ in $\Gw$.

Thus $U_{F}\rest{\Gw}$, is the unique almost large solution in
$\Gw$. \es
\Remark If $\prt \Gw$ is compact then, by
\rth{almostlarge} , $U_{\prt \Gw}$ is an almost large solution in
$\Gw$. Therefore, in this case, the last statement means that {\em
$U_{\prt \Gw}$ is the unique almost large solution in $\Gw$}.\2
\end{comment}
 \bth{un-open} Let $\Gw$ be an open bounded set in $\RN$ \sth $\Gw=\cup \Gw_n$, where $\set{ \Gw_n}$
 is an increasing family of open sets satisfying
\begin{equation}\label{un-open1}
  \C2q(\Gw\sms \Gw_n)\to 0.
\end{equation}
Put
\begin{equation}\label{un-notation}
\BAL &F_n:=\prt\Gw_n, \q D_n=\RN\sms \bar\Gw_n, \q \Gw^n=\Gw\sms \Gw_n\\
 &F:=\prt_q \Gw=\wtl\Gw\sms \Gw,\q D:=\RN\sms \wtl \Gw
 \EAL \ee
and assume that
\begin{equation}\label{un-open2}
\C2q(F_n\sms \wtl D_n)\to 0.
\ee

Under these assumptions, $V_{\wtl D}$ is the unique $\prt_q$-large solution in
$\Gw$.
\es

The proof is based on several lemmas.
\blemma{almostlarge 1} Let $\Gw$ be a bounded open set  such that,
with the notation $F:=\prt \Gw$, $D:=\RN\sms \bar\Gw$,
\begin{equation}\label{un-open4}
\C2q(F\sms \wtl D)=0. 
\ee
Then $V_{\wtl D}$ is the unique $\prt_q$-large solution in $\Gw$.
\es
\bproof
Let $v$ be a $\prt_q$-large solution in $\Gw$. First we show that
\begin{equation}\label{uniq2}
V_{D}\leq v \txt{in $\Gw$.}
\end{equation}
If $\mu\in \CT_D$ then $\mu=\sup\set{\mu\chr{K}:\, K\sbs D,\; K\text{ compact}}$
and $u_\mu=\sup u\ind{\mu\chr{K}}$ over compact sets $K$ as above.
Therefore it is sufficient to show that
\begin{equation}\label{uniq3}
u_\mu\leq v
\end{equation}
 for every $\mu\in W^{-2,q}_+(\RN)$ supported in a compact set $K\sbs D$. Since $K\cap \bar\Gw=\ems$,
$u_\mu$ is uniformly bounded in $\bar\Gw$.

Let
$$A_v:=\set{y\in  F:\,\liminf_{\Gw \ni x\to y}v(x)<\infty}.$$
Note that
$$\BAL \prt_q D&\sbs \prt D\sbs \prt\Gw=F,\\
\prt_q D&\sbs \prt_q(\RN\sms \wtl \Gw)=\prt_q\Gw.
\EAL$$
By \req{un-open4} $\C2q(F\sms \prt_q D)=0$; therefore $\C2q(F\sms\prt_q\Gw)=0$.
Therefore any $\prt_q$-large solution in $\Gw$ is an almost large solution in $\Gw$. Hence $\C2q(A_v)=0$.

Let $G_\ge$ be an open \ngh of $A_v(F)$
\sth $\C2q(G_\ge)<\ge$. Put
$$\Gw_{\gd}=\set{x\in \Gw: \dist(x,F)<\gd}, \q   \Gw'_{\gd}=\set{x\in \Gw: \dist(x,F)>\gd}.$$
Let $\Gw^*_\gd$ be a smooth domain \sth $\Gw'_\gd\sbs \Gw^*_\gd \Sbs \Gw'_{\gd/2}$. Put
$$G_{\ge,\gd}:=G_\ge\cap (\Gw\sms\ovl {\Gw^*_\gd}).$$
Then $v+V_{G_{\ge,\gd}}$ is a supersolution of \req{eq} in  $\Gw^*_\gd$ and, if $\gd$ is sufficiently small,
$$u_\mu\leq v+V_{G_{\ge,\gd}} \txon{\prt\Gw^*_\gd}.$$
Thus
$$u_\mu\leq v+V_{G_{\ge,\gd}} \txin{\Gw^*_\gd}.$$
Since $\C2q(G_\ge)\to 0$ as $\ge\to 0$,  \rth{VFest} implies that, for fixed $\gd>0$,
$$ \lim_{\ge\to0}V_{ G_{\ge,\gd}}=0 \q\text{in }\Gw^*_\gd.$$
Letting $\gd\to0$ we obtain \req{uniq3} and hence \req{uniq2}.
Further, by \rth{VtlE},
\begin{equation}\label{uniq4}
V _{\wtl D}=V_D\leq v\txt{in}\Gw.
\end{equation}

Next we show that the opposite inequality,
\begin{equation}\label{uniq4'}
 v\leq V _{\wtl D},
\ee
 is also valid.
(A-priori this is not obvious because we do not assume that $v$ is \gsmod.)

By \req{un-open4}
 $\C2q(\ovl D\sms \wtl D)=0$;
hence $V_{\ovl D}= V _{\wtl D}$.

 Let $R$ be sufficiently large so that $\ovl \Gw\sbs B_R(0)$. Then
$$(i)\q V_{\ovl D}\leq V_{\overline {D\cap B}_r}+V_{ B_R^c},\q (ii)\q U_{\ovl D}\leq U_{\overline {D\cap B}_R}+U_{ B_R^c}.$$
and $V_{\ovl D}$ (resp. $U_{\bar D}$) is the largest solution in $\Gw$,  dominated by the right hand side of inequality (i) (resp. (ii)).
Since
$\overline {D\cap B}_R$ is compact,
 $V_{\overline {D\cap B}_R}=U_{\overline {D\cap B}_R}$. The uniqueness of large solutions in smooth domains implies that
 $$U_{ B_R^c}=U_{\prt B_R}=V_{\prt B_R}=V_{B_R^c}.$$
 Combining these facts we conclude that
 $$U_{\ovl D}=V_{\ovl D}= V _{\wtl D}.$$

By definition, $v\leq U_{\ovl D}$ in $\Gw$; hence $v\leq V _{\wtl D}$.
\eproof

\blemma{almostlarge 3} Let $v$ be a solution of \req{eq} in a bounded open set $\Gw$. Suppose that $A$ is
a \qfine closed subset of $\prt\Gw$ \sth
\begin{equation}\label{partial blowup}
 \lim_{x\to y}v(x)=\infty \forevery y\in \prt \Gw\sms A.
\end{equation}
If $D:=\RN\sms \wtl\Gw$ then,
\begin{equation}\label{VD<v}
    V _{\wtl D}= V_D\leq v+ V_A\txin{\Gw}.
\end{equation}
\es
\bproof Let $\mu$ be  a measure in $W^{-2,q}_+(\RN)$ concentrated on a compact set $K\sbs D$. Let $\set{O_n}$
be a decreasing \seq of open sets \sth
$$ A\sbs O_n,\q \C2q(O_n\sms A)\to 0,\q \lim_{\Gw\ni x\to \prt \Gw\sms O_n}v(x)=\infty.$$
Let $\Gw_{n,\gd}^*$ be as in the proof of  \rlemma{almostlarge 1} and
let $\set{\gd_n}$ be a \seq of positive numbers decreasing to zero. Denote
$$G^n:=O_n\cap (\Gw_n\sms\ovl{\Gw_{n,\gd_n}^*}).$$
As in the proof of \rlemma{almostlarge 1}, we obtain
$$u_\mu\leq V_{G^n}+ v\txin{\Gw_{n,\gd_n}^*}.$$
Since $A\sbs G^n$ and $\C2q(G^n\sms A)\to 0$ it follows that $V_{G^n}\dar V_{A}$. Letting $n\tin$ we obtain
$$u_\mu\leq V_A +v$$
which in turn implies \req{VD<v}.
\eproof
\blemma{almostlarge 2} Put
$$S_{n,1}:=\ovl D_n\sms \wtl D_n,\q S_{n,2}=(\prt_q D_n)\Gd \prt_q (\RN\sms \wtl\Gw_n),\q  E_n:=F_n\Gd F. $$
Then, under the assumptions of the theorem,
\begin{equation}\label{un-open5}
\text{(a) }\; \C2q(S_{n,1})\to 0,\q \text{(b) }\;\C2q(S_{n,2})\to 0,\q \text{(c) }\;\C2q(E_n)\to 0.
\end{equation}
\es
\bproof
Since $S_{n,1}\sbs F_n\sms \wtl D_n$, (a) follows from \req{un-open2}.

Since $D_n\sbs \RN\sms \wtl\Gw_n$ it follows that
$$\prt_q D_n\sbs \prt_q (\RN\sms \wtl\Gw_n)\cup \prt_q(\bar D_n\sms \wtl D_n),\q
\prt_q (\RN\sms \wtl\Gw_n)\sbs \prt_q D_n\cup \prt_q (\bar D_n\sms \wtl D_n).$$
But \req{un-open5} (a) implies that $\C2q(\prt_q (\bar D_n\sms \wtl D_n))\to 0$. Therefore, the previous relations
imply \req{un-open5} (b).

In order to establish (c) we observe that,
$$F\sbs \prt_q\Gw_n\cup \prt_q \Gw^n,\q \prt_q \Gw_n \sbs \prt_q \Gw^n\cup F.$$
It is known that (see \cite{AH}) there exists a constant $c(N,q)$ \sth, for every Borel set $A$,
\begin{equation}\label{C(tlA)}
 \C2q(\wtl A)\leq c\C2q(A).
\end{equation}
Therefore \req{un-open1} implies that $\C2q(\wtl \Gw^n)\to0$, which in turn implies  that $\C2q(\prt_q \Gw^n)\to0$. We conclude that
\begin{equation}\label{un-open5'}
  \C2q(F\Gd \prt_q\Gw_n)\to 0.
\end{equation}
Hence, as $\prt_q \Gw_n\sbs F_n$,
\begin{equation}\label{part c1}
\C2q(F\sms F_n)\leq \C2q(F\sms\prt_q\Gw_n)\to 0.
\ee

 On the other hand,
\begin{equation}\label{un-open5*}
F_n\sms F\sbs (F_n\sms\prt_q\Gw_n)\cup (\prt_q\Gw_n\sms F). 
\ee
Since
$$\prt_q \Gw_n\supseteq\prt_q \wtl\Gw_n=\prt_q (\RN\sms \wtl\Gw_n).$$
\req{un-open5} (b) implies
$$\C2q(\prt_q D_n\sms \prt_q \Gw_n)\to 0.$$
This fact and assumption \req{un-open2} imply
\begin{equation}\label{un-temp1}
  \C2q(F_n\sms \prt_q \Gw_n)\to 0.
\end{equation}
Finally, \req{un-open5'}, \req{un-open5*} and \req{un-temp1} imply
\begin{equation}\label{part c2}
\C2q(F_n\sms F)\to 0.
\ee
This together with \req{part c1} yields \req{un-open5} (c).
\eproof
\bcom
{\em Keep in reserve}
Let $\mu$ be as before: a measure in $W^{-2,q}_+(\RN)$ concentrated on a compact set $K\sbs D$.
Since $F_n\sbs F\cup E_n$, it follows from the definition of $A_v$ (see \req{AvF}), that
$A_v(F_n)\sbs E_n\cup A_v(F)$. Since, by assumption $\C2q(A_v(F))=0$, \req{un-open5} implies that $\C2q(A_v(F_n))\to 0$.
Therefore, by \req{ZqE},
$$\ga_n:=\C2q(A_v(F_n)\cup Z^q(F_n))\to 0.$$
\Consy there exists a \seq of open sets $\set{O_n}$ \sth
$$A_v(F_n)\cup Z^q(F_n)\sbs O_n,\q \C2q(O_n)<2\ga_n.$$
Let $\Gw_{n,\gd}^*$ be a smooth domain \sth $(\Gw_n)'_\gd\sbs \Gw_{n,\gd}^* \Sbs (\Gw_n)'_{\gd/2}$. Further, let $\set{\gd_n}$
be a \seq of positive numbers decreasing to zero and denote
$$G^n:=O_n\cap (\Gw_n\sms\ovl{\Gw_{n,\gd_n}^*}).$$
\end{comment}
\vskip 3mm
\noindent{\em Proof of \rth{un-open}.}\hskip 2mm Let $A_n=F_n\sms F$. By \req{un-open5} (b), $\C2q(A_n)\to 0$.
If $v$ is a $\prt_q$-large solution in $\Gw$ then $v$ blows up \qae on $F$ and \consy it blows up \qae on $F_n\sms A_n$.
Applying \rlemma{almostlarge 3} to $v$ in $\Gw_n$ we obtain
$$   V _{\wtl D_n}=V_{D_n}\leq v+ V_{A_n}\txin{\Gw_n}. $$
Note that
$$D_n\sms D=\wtl \Gw\sms \ovl\Gw_n=(\Gw\sms \ovl\Gw_n)\cup(F\sms\ovl\Gw_n)\sbs (\Gw\sms \ovl\Gw_n)\cup(F\sms F_n)$$
and
$$D\sms D_n= \ovl\Gw_n\sms \wtl \Gw=(\Gw_n\sms \wtl \Gw)\cup(F_n\sms\wtl \Gw)\sbs (F_n\sms F).$$
Therefore, \req{un-open1} and \req{un-open5}~(c) imply that
\begin{equation}\label{CDntoCD}
  \C2q(D_n\Gd D)\to0.
\end{equation}
The definition of $V_E$ (see \req{VE}) implies
$$V _{D_n}\leq V_{D}+ V_{D_n\sms  D}$$
and, by \req{CDntoCD} and \rth{VFest}, $V_{D_n\sms  D}\to 0$. Hence
\begin{equation}\label{VDntoVD}
 V _{\wtl D_n}\to V_{\wtl D}.
 \ee
By \rth{VFest}, $V_{A_n}\to0$ in $\Gw$.
Therefore, letting $n\tin$, we obtain
$$V _{\wtl D}\leq v \txt{in $\Gw$.}$$

It remains to show that $v\leq V_{\wtl D}$. As $U_{\ovl D_n}$ is the maximal solution in $\Gw_n$,
$$ v\leq V_{\ovl D_n}=U_{\ovl D_n}.
$$
\rlemma{almostlarge 2}, implies that
$$V_{\ovl D_n}-V _{\wtl D_n}\to 0 \txt{in $\Gw$.}$$
Indeed, as an immediate consequence of the definition of $V_E$ (see \req{VE}),
$$V_{\ovl D_n}\leq V _{\wtl D_n}+ V_{\ovl D_n\sms\wtl D_n}.$$
By \req{un-open5}~(a), $\C2q(\ovl D_n\sms\wtl D_n)\to 0$. Hence, by \rth{VFest}, $V_{\ovl D_n\sms\wtl D_n}\to 0$ in
$\Gw_n$.
It follows that
$$\lim V_{\ovl D_n}\leq \lim V _{\wtl D_n}.$$
The limits exist because of monotonicity. Since $V _{\wtl D_n}\leq V_{\ovl D_n}$ we obtain,
$$\lim V_{\ovl D_n}= \lim V _{\wtl D_n}.$$
 Therefore
$$v\leq \lim V _{\wtl D_n}=V_{\wtl D}.$$
\qed

\bcom 
 =======================

Note that,
\begin{equation}\label{VDn}
\ovl D_n = \wtl D_n\cup (F_n\sms \wtl D_n) \Lra V_{\ovl D_n}\leq V_{\wtl D_n}+V_{F_n\sms \wtl D_n}.
\end{equation}
We claim that,
\begin{equation}\label{V(F-D)}
     \lim_{n\tin}V_{F_n\sms \wtl D_n}(x)=0 \forevery x\in \Gw.
\end{equation}
Indeed, given $x_0\in \Gw$, there exists $m\in \BBN$ \sth $x_0\in\Gw_m$. Choose $0<r<R$ \sth $\ovl B_r(x_0)\sbs \Gw_m$ and
$\wtl \Gw\sbs B_R(0)$ and denote $Q:=\ovl B_R(0)\sms B_r(x_0)$. Then $F_n\sbs Q$ for all $n\geq m$. Obviously $W_Q(x)<\infty$
for every $x\in B_r(x_0)$. Therefore, by \req{un-open2} and \rth{VFest},
$V_{F_n\sms \wtl D_n}(x)\to 0$ for every $x\in B_r(x_0)$. Thus \req{V(F-D)} is valid
and, by  \req{VDn},
\begin{equation}\label{VDn1}
0\leq (V_{\ovl D_n}- V_{\wtl D_n})(x)\to 0 \q\text{as }n\tin, \forevery x\in \Gw.
\end{equation}

=========================

As $v_n\leq V_{\ovl D_n}$ in $\Gw_n$ we obtain
\begin{equation}\label{tlDn}
 \limsup_{n\tin} (v_n-V_{\wtl D_n})\leq 0 \txin{\Gw}.
\end{equation}
 The function $v^*_n:=(v+V\ind{F_n\cap\Gw})\rest{\Gw_n}$
 is a supersolution in $\Gw_n$ which  dominates $v'_n$ and \consy dominates $v_n$. Hence,
\begin{equation}\label{VDn2}
 v\leq v_n\leq v_n + V\ind{F_n\cap\Gw} \txin{\Gw_n}.
\end{equation}
Note that
 $\Gw\sbs\RN\sms \wtl D$ and \consy
 $$F_n\cap\Gw\sbs F_n\sms \wtl D\sbs F_n\sms F.$$
By \req{un-open5} and \rth{VFest},
$$V\ind{F_n\cap\Gw}\leq V_{F_n\sms F}\to 0.$$
 Therefore, by \req{VDn2},
\begin{equation}\label{VDn3}
 \lim v_n=v \txin{\Gw}
\end{equation}
and further, by \req{tlDn},
\begin{equation}\label{VDn4}
 v\leq \lim V_{\wtl D_n} \txin{\Gw}.
\end{equation}
The limit above exists because $\set{D_n}$ is decreasing. In addition,
$$D_n\sms D= \wtl \Gw\sms \ovl \Gw_n = (F\cup \Gw)\sms(F_n\cup\Gw_n)\sbs (F\sms F_n)\cup(\Gw\sms \Gw_n)$$
and, by \req{un-open1} and \req{un-open5},
$$\C2q(D_n\sms D)\to 0 \txt{and} V_{D_n}\dar V_D.$$
In view of \rth{VtlE}, this is equivalent to
\begin{equation}\label{VDn5}
 V_{\wtl D_n}\dar V_{\wtl D} \txin{\Gw}.
\end{equation}
Finally, by \req{VDn4}, \req{VDn5} and \req{VD<v} we obtain $v= V_{\wtl D}$.
\end{comment}
\bcom
\bcor{uniq1} Let $\Gw$ be a bounded domain and put $F=\prt \Gw$ and $D=\RN\sms \ovl \Gw$.
If \qae point  $y\in F$ satisfies the
Wiener criterion
\begin{equation}\label{WFD1}
   W_{D}(y)=\infty,
\end{equation}
there exists a unique almost large solution in $\Gw$. In addition
$V_{\wtl D}=V_{\ovl D}=U_F$ in $\Gw$.

If {\em every} point  $y\in F$ satisfies \req{WFD1}
then $V_{\ovl D}$ is a large solution in $\Gw$ and
\begin{equation}\label{WFD2}
   W_{F}(y)=\infty \forevery y\in F.
\end{equation}
\es
\note{Remark} Condition \req{WFD2} holds in every Lipschitz domain $\Gw$.
\bproof The definition of $W_D$ implies that $W_D\sim W_{\wtl D}$. Therefore \req{WFD1} implies that
$W_{\wtl D}(y)=\infty$ for $\C2q$ a.e. $y\in F$.
By \rth{VFest}, this implies that
$V_{\wtl D}$ blows up $\C2q$ a.e. at $F$. We claim that
\begin{equation}\label{WFD3}
 \C2q(F\sms {\wtl D})=0.
\end{equation}
By negation, suppose that $\C2q(F\sms {\wtl D})>0$ and let $K\sbs F\sms {\wtl D}$ be a compact set of positive $\C2q$ capacity
\sth $W_D(y)=\infty$ for every $y\in K$. Then, by \rth{VtlE}, $(U_K-V_{\wtl D})_+\not\equiv 0$ so that $U_K+V_{\wtl D}>V_{\wtl D}$.
On the other hand $U_K+V_{\wtl D}\leq U_{\ovl D}$.

By \rth{un-open}, $V_{\wtl D}$ is the unique almost large solution in $\Gw$. Since $U_F\rest{\Gw}$ and $U_{\ovl D}$ are also
almost large solutions in $\Gw$ it follows that $V_{\wtl D}=V_{\ovl D}=U_F$ in $\Gw$.
solutions of \req{eq} in $\Gw$ which blow up $\C2q$ a.e. at $F$.

In addition \req{WFD1} implies that
$\C2q(\ovl D\sms \wtl D)=0$.
Therefore the first assertion is an immediate \cons of \rth{un-open} and \rth{low-est}.
The fact that $V_{\ovl D}=U_F\rest{\Gw}$
is a \cons of the uniqueness of the almost large solution. The last assertion of the corollary follows from this fact
and the Theorems \ref{t:low-est} and \ref{t:up-est}.
\eproof
\end{comment}
\bcor{uniq2} Suppose that $\Gw=\cup_1^\infty Q_n$ where $\set{Q_n}$ is a \seq of  open sets \sth
\begin{equation}\label{cor-un1}
\sum_1^\infty\C2q(Q_n)<\infty.
\ee
For every $n\in\BBN$, put
$$S_n=\cup_1^n \ovl Q_k,\q D_n=\RN\sms \ovl S_n$$
and assume that
\begin{equation}\label{WFD'}
  \C2q(\prt S_n\sms \wtl D_n)\to 0.
\end{equation}
Then there exists a unique almost large solution in $\Gw$.
\es
\note{Remark} If $y\in \prt S_n$ and there exists an open cone $C_y$, with vertex $y$, \sth $C_y\sbs \RN\sms S_n$ then
$y\in \prt_q S_n$. Hence if, for every $n\in \BBN$, this condition is satisfied \qae on $\prt S_n$ then \req{WFD'} holds.
In particular, if $\set{Q_n}$ is a \seq of balls, \req{WFD'} is satisfied.
\bproof Let $\Gw_n= S_n^0:=S_n\sms \prt S_n$. Then $\set{\Gw_n}$ is an increasing \seq of open sets, \req{WFD'} implies
\req{un-open2} and \req{cor-un1} implies \req{un-open1}. Therefore
the corollary is an immediate consequence of \rth{un-open}.
\eproof
\bcom
By \req{WFD'}, $\prt Q_n\sbsq \wtl E_n=\prt_q E_n\cup E_n$. As $E_n:=\BBR^N\sms \ovl Q_n$ we obtain
 $\prt Q_n\sbsq \prt_q E_n$ and $\prt_q E_n=\prt_q\ovl Q_n$.
 Therefore
$$\prt\Gw_n\sbs  \cup_1^n \prt Q_k\sms S_n^0 \sbsq  \cup_1^n \prt_q\ovl Q_n\sms S_n^0
\sbsq \prt_q S_n\cup (\qinterior S_n\sms S_n^0).$$
But $ \qinterior S_n\sms S_n^0=\prt S_n\sms \prt_q S_n$ and, if $x\in \prt S_n\sms \prt_q S_n$ then $x\in \prt Q_j$ for some $j$ but
$x\not\in \prt_qS_n^c=\prt_q(\cap_1^n E_k)$.
Thus \req{un-open2} is valid in the present case. In addition \req{cor-un1} implies \req{un-open1}. Therefore
the corollary is a consequence of \rth{un-open}.
\eproof
\end{comment}
\note{Example} Let $\set{x^m}$ be a \seq of distinct points in $B_1(0)$. Let $\set{r_n}$ be a decreasing \seq of positive numbers
\sth $\set{B_{r_n}(x^n)}$ is a \seq of  balls contained in $B_1(0)$ and
\begin{equation}\label{cor-un2}\BAL
  &\sum r_n^{N-2q'}<\infty &&\text{if }N>2q',\\
 &\sum (1-\log r_n)^{1-q'}<\infty &&\text{if }N=2q'.\EAL
\end{equation}
Then there exists a unique large solution in $\Gw:=\cup_1^\infty B_{r_n}(x^n)$.

Indeed $\C2q(B_r)\sim r^{N-2q'}$ if $N>2q'$ and $\C2q(B_r)\sim \log (1-\log r)$ if $N=2q'$ and $0<r<1$. Therefore the conditions
of \rcor{uniq2} are satisfied.

Note that $\wtl\Gw=\cup\ovl B_{r_n}(x^n)$, but, in general $\ovl \Gw$ is much larger. For instance,
if $\set{x^m}$ is a dense \seq in $B_1(0)$ then $\ovl \Gw=\ovl B_1(0)$.
 Therefore it
is important that our conditions in \rcor{uniq2} require $\C2q(\prt_q \Gw\sms\wtl D)=0$ and not $\C2q(\prt \Gw\sms\wtl D)=0$.
\section{Very weak subsolutions}\label{sec:weaksol}
In this section $F$ is a \qfine closed  set contained in $B_1(0)$ and $D=B_2(0)\sms F$. Note that $D$ is a \qfine open set, but
not necessarily open in the Euclidean topology.

We denote by $W^{2,q'}(D)$ the set $\set{h\rest{D}:h\in \W2q(\RN)}$. If $f\in W^{2,q'}(\RN)$ we denote by $\qsupp f$ $(=$
the $\C2q$-fine support of $f)$ the
intersection of all \qfine closed sets $E$ \sth $f=0$ a.e. in $\RN\sms E$.

The following subspace of $W^{2,q'}(D)$ serves as a space of test
functions  in our study:
\begin{equation}\label{testfunction}
W_{0,\infty}^{2,q'}(D):= \set{h\rest{D}: \,h\in \W2q(\RN)\cap L^\infty(\RN),\; \qsupp h\Sbs D}.
\end{equation}
The notation $E\Sbs D$ means: $E$ is `strongly contained' in
$D$, i.e., $\ovl E$ is a compact subset of $D$. Some features of this space are discussed in Appendix~A.


The following statement was established in \cite{MV3} (see Lemma 2.6). (The framework in \cite{MV3}
is somewhat different, but the proof, with obvious modifications, applies to the present case as well.)
\blemma{qexhaust}
Let $D$ be a bounded \qfine open set. Then there exists an increasing \seq of compact sets $\set{E_n}$ \sth
\begin{equation}\label{qexhaust}\BAL
 & E_n\sbs \mathrm{int}_q E_{n+1},&&\cup_1^\infty E_n\sbs D,\\
 & \C2q(D\sms \cup_1^\infty E_n)=0,&&\C2q(E_n)\to \C2q(D).
\EAL\end{equation}
\es

A  \seq  of sets $\set{E_n}$ as above is called a {\em $q$-exhaustion of $D$}.

\bcom
Moreover, there exists a sequence of functions $\set{\gf_n}$ in
$W_{0,\infty}^{2,q'}(D)$ \sth
\be\label{basicseq}
 0\leq \gf_n\leq 1, \q \gf_n \leq \gf_{n+1},\;n=1,2,\ldots,\q \cup_n \qsupp\gf_n=D.
\ee
Thus the first statement holds with $E_n=\qsupp\gf_n$.
\end{comment}
\vskip 1mm

We denote by $L^q\qloc(D)$  the space of measurable functions $f$ in
$D$ \sth, for every positive $\phi\in W_{0,\infty}^{2,q'}(D)$, $f\in
L^q(D;\phi)$, i.e., $|f|^q\phi\in L^1(D)$. We endow this space with
the topology determined by the family of semi-norms
\begin{equation}\label{seminorms}
  \set{\norm{\cdot}\ind{L^q(D;\phi)}:\,  \phi\in
W_{0,\infty}^{2,q'}(D),\;\phi\geq 0}.
\end{equation}
This topology will be denoted by $\tau\qloc(D)$.

Further we denote by $\GTM\qloc(D)$  the space of positive Borel measures $\mu$ in
$D$ \sth
\begin{equation}\label{MUq}
\begin{aligned}
(a)\qq&K\sbs D, \;\text{$K$ compact}&&\Lra \mu(K)<\infty, \\
 (b)\qq&E\sbs D,\; E \text{ Borel, }\C2q(E)=0&&\Lra \mu(E)=0.\\
 \end{aligned}
\end{equation}

We observe that,
\blemma{int1}
 If $\mu\in\GTM\qloc(D)$ then:\\
 (i) There exists an increasing \seq $\set{\mu_n}$ of positive,
bounded measures in $W^{-2,q}(\RN)$ \sth $\mu_n(D^c)=0$ and $\mu_n\uparrow \mu$.\\
(ii) $W_{0,\infty}^{2,q'}(D)\sbs L^1(\mu)$.
\es
\bproof
(i) This is well known in the case that $\mu$ is a positive, bounded measure \cite{BP} and it
 follows from \rlemma{qexhaust}
in the case that $\mu$ is a positive measure in $\GTM\qloc(D)$.\2
(ii)  If $\vgf\in W_{0,\infty}^{2,q'}(D)$,
it vanishes outside a compact set $K_\vgf\sbs D$. By definition, $\mu(K_\vgf)<\infty$.
Furthermore $\vgf$ is the limit \qae of smooth functions; \consy it is $\mu$-measurable.
Since $\vgf$ is bounded, it is integrable relative to $\mu$.
\eproof
 \note{Notation} A \seq $\set{\mu_n}$  as in \rlemma{int1}~(i)
will be called  a {\em determining \seq} for $\mu$.

We introduce below a very weak type of
subsolution of \req{eqmu} defined as follows.
\bdef{vwsubsol} {\rm Assume that the measure $\mu$ in \req{eqmu} belongs to $\GTM\qloc(D)$.
A non-negative measurable function $u$ is a {\em very weak subsolution} of \req{eqmu} in $D$ if, for every
non-negative $\gf\in W^{2,q'}_{0,\infty}(D)$,}
\begin{align}\label{vwsub1}
 & u \in L^q(D;\gz)\txt{where} \gz:=\gf^{2q'},\\
-&\int_D u\Gd \gz dx + \int_D u^q\gz dx\leq\int_D\gz d\mu.\label{vwsub2}
\end{align}
\es
\note{Remarks} (a) If \req{vwsub1} holds for every non-negative $\gf\in W^{2,q'}_{0,\infty}(D)$ then
\begin{equation}\label{vwsub1'}
u\Gd \gz\in L^1(D).
 \ee
This is proved in the next lemma.\1
 (b) Let $\gf\in W^{2,\gg}_{0,\infty}(D)$,  $\gg\geq 1$.
 By interpolation, $|\nabla\gf|^2\in L^{\gg}(D)$ and
\begin{equation}\label{gradphi}
 \norm{|\nabla \gf|^2}\ind{L^{\gg}(D)}\leq c(q,N)L\norm{D^2\gf}\ind{L^{\gg}(D)},\q L:=\norm{\gf}\ind{L^\infty(D)}.
\end{equation}
where $|D^2\gf|:=\sum_{|\ga|=2}|D^\ga\gf|$.\1
(c) If $\gf\in W^{2,\gg}_{0,\infty}(D)$,  $\gg\geq 1$ then
\begin{equation}\label{Gdphi}\BAL
  &|\gf|^{2\gg}=(\gf^{2})^{\gg}\in \W2q_{0,\infty}(D),\\
  &\nabla (|\gf|^{2\gg})=2\gg(\gf^{2})^{\gg-1/2}\nabla\gf,\\
  & \Gd(\gf^{2\gg})=2\gg(2\gg-1)|\gf|^{2\gg-2}|\nabla\gf|^2+2\gg|\gf|^{2\gg-1}\Gd\gf.
\EAL\end{equation}
The last two formulas are easily verified for $\gf\in C_c^\infty(\RN)$; in the general case they are obtained
 by the usual density argument. The fact that $|\gf|^{2\gg}\in \W2q_{0,\infty}(D)$ is a \cons of these
 formulas and \req{gradphi}. \req{Gdphi} imply,
 \begin{equation}\label{normphi}
 \norm{|\gf|^{2\gg}}\ind{W^{2,\gg}(D)}\leq AL^{2\gg-1}\max(1,L) \norm{\gf}\ind{W^{2,\gg}(D)}.
\end{equation}
\bth{int<W2qnorm}
$\mathrm{(i)}$ If $u$ is a non-negative measurable function satisfying \req{vwsub1}
then $u\Gd\gz\in L^1(D)$ .\1
$\mathrm{(ii)}$ If
 $u$ is a very weak subsolution of \req{eq} in $D$ (i.e. $\mu=0$) then, for every non-negative $\gf\in \W2q_{0,\infty}(D)$,
\begin{equation}\label{W2qnorm0}
 \int_D u|\Gd\gz|dx + \int_D u^q\gz\leq c
 \left(L\norm{D^2\gf}\ind{L^{\gg}(D)}\right)^{\gg},
\end{equation}
where $\gz:=\gf^{2q'}$, $c=c(N,q)$ and $ L:=\norm{\gf}\ind{L^\infty(D)}$.\1
$\mathrm{(iii)}$ Let $\mu\in W^{-2,q}(\RN)$ be a positive bounded
measure vanishing outside $D$.
If $u$ is a non-negative very weak subsolution of \req{eqmu} then
\begin{equation}\label{W2qnorm1}\BAL
 \norm{u}\ind{L^q(D,\gz)}\leq
 &cL^{1/(q-1)}\Big(\big(\norm{D^2\gf}\ind{L^{q'}(D)}\big)^{1/(q-1)} +\\
&\big(L^{1/(q-1)}\norm{\mu}\ind{W^{-2,q'}}\norm{D^2\gf}\ind{L^{q'}(D)}\big)^{1/q}\Big)
\EAL
\end{equation}
Finally, if $L\leq \ovl L$,  $u$ satisfies
\begin{equation}\label{W2qnorm1'}\BAL
 &\norm{u}\ind{L^q(D,\gf)}\leq\\
 &c(N,q,\ovl L)L^{\rec{q-1}}\Big(\norm{D^2\gf}^{\rec{q-1}}\ind{L^{q'}(D)} +
\norm{\mu}^{\rec{q}}\ind{W^{-2,q'}}\norm{D^2\gf}^{\rec{q}}\ind{L^{q'}(D)}\Big).
\EAL\end{equation}
\es
\bproof Let $\gf$ be a non-negative function in $W^{2,q'}_{0,\infty}(D)$. By  \req{Gdphi}, with $\gg=q'$, we obtain
\bcom
\begin{equation}\label{W2qnorm2}
  \BAL \Gd\gz&=2q'\gf^{2q'-1}\Gd\gf+2q'(2q'-1)\gf^{2q'-2}{\absol{\nabla\gf}}^{2}\\
&=\gf^{2q'-2}\big(2q'(2q'-1){\absol{\nabla\gf}}^{2}+2q'\gf\Gd\gf\big),
\EAL
\end{equation}
so that
\end{comment}
$$\absol{\Gd\gz}\leq c(q) \gz^{1/q}M(\gf),\q M(\gf):=\big(\absol{\nabla\gf}^{2}+\gf\absol{\Gd\gf}\big)$$
and hence, using \req{gradphi},
\begin{align}\label{W2qnorm3}
&\int_{D}u\absol{\Gd\gz} dx
\leq c(q)\left(\int_{D}u^q\gz\,dx\right)^{1/q}\left(\int_{D}M(\gf)^{q'}dx\right)^{1/q'},\\
&\int_D M(\gf)^{q'}dx \leq c(q,N)L\norm{D^2\gf}^{q'}\ind{L^{q'}(D)}. \label{W2qnorm4}
\end{align}
Assuming that $u\in L^q(D,\gz)$ we obtain $u\Gd\gz\in L^1(D)$.
\bcom
Now let $\psi$ be a non-negative function in $\w02q$ and put
$$\gf:=(1+\psi)^{\rec{2q'}}-1,\q \gz:=\gf^{2q'}.$$
 Then  $\gf\in \w02q$ and $u\Gd\gz\in L^1(D)$. Furthermore, $\gz\sim \psi$ and
 \begin{equation}\label{vwsub1''}
\norm{D^2\gf}\ind{L^{q'}(D)}\leq c(N,q)\norm{D^2\psi}\ind{L^{q'}(D)}(1+\norm{\psi}\ind{L^{\infty}(D)}).
 \end{equation}
Hence, by \req{Gdphi},
$$\Gd(\gf^{2\gg})=2\gg(2\gg-1)|\gf|^{2\gg-2}|\nabla\gf|^2+2\gg|\gf|^{2\gg-1}\Gd\gf.$$
$$\Gd\gz=2\gg(2\gg-1)|\gf|^{2\gg-2}|\nabla\gf|^2+2\gg|\gf|^{2\gg-1}\Gd\gf.$$
$$2\gg|\gf|^{2\gg-1}\Gd\gf= \Gd\gz-2\gg(2\gg-1)|\gf|^{2\gg-2}|\nabla\gf|^2.$$

 it follows that $u\Gd\psi\in L^1(D)$.

\end{comment}
\vskip 1mm
We turn to the proof of (ii) and (iii).  We assume that $u$ is a non-negative very weak subsolution as in \rdef{vwsubsol}.
Put
$$A:=\left(\int_Du^q\gz dx\right)^{1/q},\; B:=\left(\int_D M(\gf)^{q'} dx\right)^{1/q'},\;
C:=\norm{\mu}\ind{W^{-2,q'}}\norm{\gz}\ind{W^{2,q'}(D)}$$
By \req{vwsub2} and \req{W2qnorm3}
\begin{equation}\label{W2qnorm5}
 A^q= \int_D u^q\gz dx\leq \int_D u\Gd \gz dx +\int_D\gz\,d\mu\leq c(q,N)AB+C.
\end{equation}
This implies
$$A^q\leq \rec{q}A^q+\rec{q'}(cB)^{q'}+C \Lra A^q\leq (cB)^{q'}+q'C\leq c'(N,q)\max(B^{q'},C).$$
Thus
\begin{equation}\label{W2qnorm6}
 A\leq c(q,N)\left(B^{1/(q-1)}+ C^{1/q}\right).
\end{equation}
By Poincar\'e's inequality
$$\norm{\gz}\ind{W^{2,q'}(D)}\leq c(q,N)\norm{D^2\gz}\ind{L^{q'}}$$
and therefore, by the same computation as in \req{Gdphi},
$$\norm{\gz}\ind{W^{2,q'}(D)}\leq
c(q,N)\Big(\norm{\gf^{2/(q-1)}(\nabla\gf)^2}\ind{L^{q'}}+\norm{\gf^{(1+q)/(q-1)} D^2\gf}\ind{L^{q'}}\Big).$$
Therefore by \req{gradphi} and \req{W2qnorm4},
\begin{equation}\label{W2qnorm7}
\norm{\gz}\ind{W^{2,q'}(D)}\leq L^{\frac{1}{q}}(L+L^{\rec{q}})\norm{D^2\gf}\ind{L^{q'}}.
\ee
This estimate  and \req{W2qnorm6} imply \req{W2qnorm1}. Further, if $\mu=0$, \req{W2qnorm1}, \req{W2qnorm3} and \req{W2qnorm4}
imply \req{W2qnorm0}.

Now let $\psi$ be a non-negative function in $\w02q(D)$ and put
$$\gf:=(1+\psi)^{\rec{2q'}}-1,\q \gz:=\gf^{2q'}.$$
 Then  $\gf\in \w02q(D)$, $\gz\sim \psi$ and
 $$\norm{D^2\gf}\ind{L^{q'}(D)}\leq c(N,q)\norm{D^2\psi}\ind{L^{q'}(D)}(1+\norm{\psi}\ind{L^{\infty}(D)}).
 $$
 This inequality and \req{W2qnorm1} imply \req{W2qnorm1'}.

\eproof
\blemma{vwsub0} If $F$ is a Borel set \sth $\C2q(F)=0$ then the only non-negative very weak subsolution of \req{eq} in $D=F^c$ is
the trivial solution.
\es
\Remark A set of capacity zero is \qfine closed by defintion. Therefore the notion of very weak subsolution in $F^c$ is well defined
in the present case.
\bproof
Since $\C2q(F)=0$, there exists a \seq $\set{\eta_n}$ in $\W2q(\RN)$ \sth
$0\leq \eta_n\leq 1$, $\norm{\eta_n}\ind{\W2q}\to 0$ and $\eta_n=1$ on a \ngh of $F$ (depending on $n$).
 Applying \req{W2qnorm0} to $u$ and $\gf_n=1-\eta_n$ yields:
$$  \int_D u^q\gf_n^{2q'}dx\leq c \norm{D^2\eta_n}\ind{L^{q'}(D)}^{q'}\to 0.$$
Since $\norm{\eta_n}\ind{L^{q'}}\to 0$ it follows that there exists a \sseq converging to zero a.e.. Therefore
$$\liminf \int_D u^q\gf_n^{2q'}dx\geq \int_D u^q\,dx.$$
This implies that $u=0$ a.e.
\eproof

\section{$C_{2,q'}$-strong solutions in finely open sets and uniqueness II}\label{fineq}
We start with the definition of '$\C2q$-strong' solutions of
\req{eqmu} in a \qfine open set or more generally in a $\C2q$-quasi
open set. We recall that a set $E$ is  $\C2q$-quasi open if, for
every $\ge>0$ there exists an open set $O$ \sth $E\sbs O$ and
$\C2q(O\sms E)<\ge$. Every \qfine open set is $\C2q$-quasi open; if
$E$ is $\C2q$-quasi open then $E\qsim \qinterior E$, (see
\cite[Chapter 6]{AH}).

\bdef{finesolution} {\rm Let $D$ be a $\C2q$-quasi open set, let
$\mu\in\GTM\qloc(D)$ be a non-negative measure and let $\set{\mu_n}$
be a determining \seq for $\mu$ (see \rlemma{int1}).\1
(i) A  positive function $u\in L^q\qloc(D)$ is a {\em $\C2q$-strong
solution} of \req{eqmu} in $D$ if there exists a decreasing \seq of
open sets $\set{\Gw_n}$, \sth $D\sbs \Gw_n$
and, for each $n$,  there exists a positive solution
$u_n\in L^q\loc(\Gw_n)$ of the equation
\begin{equation}\label{un-Gwn}
  -\Gd u_n + u_n^q=\mu_n
\end{equation}
 \sth}
\begin{equation}\label{lim-Lqloc}
   u_n\to u \txin{L^q\qloc(D)}.
\end{equation}
{\rm We say that $\set{(u_n,\Gw_n)}$ is} a determining \seq for $u$
in $D$.\1 {\rm (ii) A $\C2q$-strong subsolution is defined in the
same way as above except that $u_n$ is only required to be a
subsolution of \req{un-Gwn} in $\Gw_n$.}\1 {\rm (iii) A positive
$\C2q$-strong solution of \req{eqmu} in $D$ is {\em \gsmod} if, in
addition, the \seq $\set{u_n}$ is non-decreasing and there exists a
\seq $\set{v_n}$ \sth $v_n\in L^1(\Gw_n)$ and
\begin{equation}\label{modsol}
 -\Gd v_n=\mu_n, \q u_n\leq v_n \txin{\Gw_n}, \q n=1,2,\ldots\,.
\end{equation}
{\rm (iv)} If $\mu$ is bounded and
$\set{\norm{v_n}\ind{L^1(\Gw_n)}}$ is bounded we say that $u$ is a
{\em moderate} solution. } \es \note{Remark} If $D$ is an open set
we may choose $\Gw_n=D$ for every $n$. Therefore any non-negative
solution of \req{eqmu} in $D$ is a $\C2q$-strong solution in $D$.
Furthermore, if $u$ is a \gsmod  solution of \req{eq} in $D$ in the
standard sense (i.e. the limit of an increasing \seq of moderate
solutions) then it is a \gsmod $\C2q$-strong solution in the sense
of part (iii) of the above definition.

\bdef{boundary data} {\rm (a) A $\C2q$-strong solution $v$ of
\req{eq} in $D$ is called a {\em \prtql solution}  if
\begin{equation}\label{blowup-cond}
 \qlim_{x\to \prt_q D}v(x)=\infty \q\text{\qae at}\; \prt_q D.
\ee
This condition is understood as follows. There exists a determining \seq $\set{(v_n,\Gw_n)}$ for $v$
in $D$ \sth
\begin{equation}\label{Gwn-cond}
 \sum_1^\infty \C2q(\Gw_n\sms D)<\infty,
\end{equation}
and, for every $M>0$, $k\in \BBN$, there exists an open set $Q_{k,M}$ \sth
\begin{equation}\label{QkM-cond}
\BAL  \cup_{n=k}^\infty & \wtl{\Gw_n\sms D}\sbs Q_{k,M},\q  \lim_{k\tin} \C2q(Q_{k,M})=0,\\
 &\liminf_{\BSA{c} x\to \prt_q D\sms Q_{k,M}\\x\in \Gw_n\ESA}v_n(x)\geq M \forevery n\geq k.
\EAL\ee
Note that $\prt_q D\sms Q_{k,M}\sbs \prt \Gw_n$ for all $n\geq k$.

If $F$ is a quasi closed subset of $\prt D$, the condition
\begin{equation}\label{infinity at F}
\qlim_{x\to F}v(x)= \infty \q\text{\qae at}\; F
\ee
is defined in the same way except that the second line in \req{QkM-cond} reads
$$\liminf_{x\to F\sms Q_{k,M}}v_n(x)\geq M \forevery n\geq k.$$
(b) Let $v$ be a non-negative $\C2q$-strong subsolution of \req{eq}
in $D$.
The condition
\begin{equation}\label{zero-cond}
\qlim_{x\to \prt D}v(x)= 0 \q\text{\qae at}\; \prt_q D
\ee
is understood as follows. There exists a determining \seq $\set{(v_n,\Gw_n)}$ for $v$
in $D$ satisfying \req{Gwn-cond} and a family of open sets
$$\set{Q_{k,\ge}:\, \ge>0, \; k\in \BBN}$$
 \sth,
\begin{align}\label{Qke-cond}
& \cup_{n=k}^\infty  \wtl{\Gw_n\sms D}\sbs Q_{k,\ge},\q  \lim_{k\tin} \C2q(Q_{k,\ge})=0,\\
 &\limsup_{\BSA{c} x\to \prt_q D\sms Q_{k,\ge}\\x\in \Gw_n\ESA}v_n(x)\leq\ge \forevery n\geq k.\label{Qke-0}
\end{align}

If $F$ is a quasi closed subset of $\prt D$, the condition
\begin{equation}\label{zero at F}
\qlim_{x\to F}v(x)= 0 \q\text{\qae at}\; F
\ee
is defined in the same way except that  \req{Qke-0} is replaced by
\begin{equation}\label{Qke at F}
\limsup_{x\to F\sms Q_{k,\ge}}v_n(x)\leq\ge \forevery n\geq k.
\end{equation}
(c) Let $v$ be a non-negative $\C2q$-strong solution of \req{eq} in
$D$ and let $u$ be a non-negative classical solution in a domain
$G\supseteq D$. We say that $u\qleq v$
at $\prt_q D$ 
if there exists a determining \seq $\set{(v_n,\Gw_n)}$ for $v$
in $D$ satisfying \req{Gwn-cond} and a family of open sets $\set{Q_{k,\ge}:\,\ge>0,\;k\in \BBN}$
satisfying \req{Qke-cond} \sth
\begin{equation}\label{Qke-ineq}
\limsup_{\BSA{c} x\to \prt_q D\sms Q_{k,\ge}\\ x\in \Gw_n\ESA}(u-v_n)(x)\leq\ge \forevery n\geq k.
\ee

If $F$ is a quasi closed subset of $\prt D$, the condition $u\qleq w$  at
$F$ is defined in the same way except that \req{Qke-ineq} is replaced by
\begin{equation}\label{Qke-ineq F}
  \limsup_{x\to  F\sms Q_{k,\ge}}(u-v_n)(x)\leq\ge \forevery n\geq k.
\end{equation}
}
\es

We present several results concerning $\C2q$-strong solutions. The
main ingredients in these proofs are: \rth{int<W2qnorm}, the results
of Section 5 concerning $V_F$
and the results of Appendix \ref{w0q}. 

\bth{fine-strong} For every $\ovl L>0$, there exist constants
$c=c(N,q)$ and $c=c(\ovl L)$ \sth, for every non-negative measure
$\mu\in\GTM\qloc(D)$ and every non-negative $\C2q$-strong solution
$u$ of \req{eqmu} in $D$ the following  holds:
 \begin{equation}\label{strongest}\BAL
 &\norm{u}\ind{L^q(D,\gf)}^q\leq\\
 &c(N,q)\Big((\norm{\gf}\ind{L^\infty(D)}\norm{D^2\gf})^{q'}\ind{L^{q'}(D)} + c(\ovl L)
\norm{\mu}\ind{W^{-2,q'}}\norm{D^2\gf}\ind{L^{q'}(D)}\Big)
\EAL
\end{equation}
for every $\gf\in \w02q(D)$ \sth $0\leq \gf\leq \ovl L$.

Every $\C2q$-strong solution $u$ as above satisfies,
\begin{align}\label{weaksol}
 & u^q\gf\in L^1(D),\q u\Gd (\gf\psi)\in L^1(D)\\
-&\int_D u\Gd (\gf\psi) dx + \int_D u^q(\gf\psi) dx=\int_D\gf\psi\, d\mu\label{weakeq}
\end{align}
for any non-negative $\gf,\psi\in \w02q(D)$. Finally $u$ satisfies the estimate
\begin{equation}\label{u<WF}
  u\leq c(N,q)W_F \txt{a.e. in $D$.}
\end{equation}
\es
\bproof We use the notation of \rdef{finesolution}. If $u_n$ is a solution of
\req{eqmu} in $\Gw_n$ and $\gf\in\w02q(D)$ then
\begin{equation}\label{eqD}
 \int_D u_n\Gd \gf dx + \int_D u_n^q\gf dx=\int_D\gf\, d\mu_n
\end{equation}
Evidently, $u_n$ is, in particular, a very weak subsolution in $D$; \consy it satisfies inequality \req{W2qnorm1'}.
By assumption, $u_n\to u$ in
$L^q\qloc$; hence $u$ satisfies \req{W2qnorm1'}.

Assume that  $\gf,\psi\in \w02q(D)$ and $0\leq\gf \leq \ovl L$ and the same for $\psi$.
 Clearly, \req{weakeq} holds for $u_n$. In addition,
$$\int_D u_n^q(\gf\psi) dx\to \int_D u^q(\gf\psi) dx,$$
$$\int_D(\gf\psi)\, d\mu_n\to \int_D(\gf\psi)\, d\mu.$$
Further,
$$\Gd (\gf\psi)=\gf\Gd\psi+ \psi\Gd\gf +2 \nabla\gf\cdot\nabla\psi,$$
so that
$$ \int_D u_n|\Gd (\gf\psi)| dx \leq \int_D u_n(\gf|\Gd\psi|+ \psi|\Gd\gf|)\,dx + 2\int_D u_n |\nabla\gf\cdot\nabla\psi|dx.$$
Using again the fact that $u_n\to u$ in $L^q\qloc$
$$\int_D u_n(\gf\Gd\psi+ \psi\Gd\gf)\,dx\to \int_D u(\gf\Gd\psi+ \psi\Gd\gf)\,dx.$$
In addition,
$$\BAL \int_D u_n |\nabla\gf\cdot\nabla\psi|dx&\leq (\int_D u_n(\nabla\gf)^2dx)^{1/2}(\int_D u_n(\nabla\psi)^2dx)^{1/2}
\EAL$$
By \req{eqD}
 $$\int_D u_n\Gd \gf^2 dx + \int_D u_n^q\gf^2 dx=\int_D\gf^2\, d\mu_n$$
so that
$$\int_D u_n(\nabla \gf)^2 dx \leq \rec{2}\int_D\gf^2\, d\mu_n + \int_D u_n\gf|\Gd\gf|\,dx.$$
By Fatou, this implies,
$$\int_D u(\nabla \gf)^2 dx \leq \rec{2}\int_D\gf^2\, d\mu + \int_D u\gf|\Gd\gf|\,dx.$$

Now assume temporarily that $\set{u_n}$ is non-decreasing so that $u_n\uar u$. Then, by the dominated convergence
theorem,
\begin{equation}\label{unabla}
   \int_D u_n (\nabla\gf\cdot\nabla\psi)dx \to \int_D u (\nabla\gf\cdot\nabla\psi)dx.
\end{equation}
The convergence results obtained above and \req{eqD} imply \req{weaksol} and \req{weakeq}. This in turn implies,
by \rth{int<W2qnorm}, the estimate \req{strongest}.

Discarding the assumption of monotonicity, put $v_n:=\max(u_1,\ldots,u_n)$. Then $v_n$ is a subsolution of the equation
$$-\Gd v+v^q=\mu_n \txin{\Gw_n}$$
and there exists a  solution $\bar v_n$ of this equation which is the smallest among those dominating $v_n$.
Then $\set{\bar v_n}$ is non-decreasing and, by \rth{int<W2qnorm},
$$\sup_n \int_D \bar v_n^q\gf\,dx<\infty$$
for any non-negative $\gf\in \W2q_{0,\infty}(D)$. Therefore $w=\lim
\bar v_n\in L^q\qloc(D)$ and by the previous part of the proof $w$
is a $\C2q$-strong solution in $D$. In particular
$$w|\nabla \gf||\nabla \psi|\in L^1(D) \forevery \gf,\psi\in \W2q_{0,\infty}(D).$$
Clearly,
$$u_n|\nabla \gf\cdot\nabla \psi|\leq w|\nabla \gf||\nabla \psi|.$$
Therefore, once again by the dominated convergence theorem, we
obtain \req{unabla} which together with the previous convergence
results imply \req{strongest}, \req{weaksol} and \req{weakeq}.

Put $F_n=B_R\sms\wtl \Gw_n$, $F=B_R\sms \wtl D$. In order to prove
the last assertion, we observe that, in $\Gw_n$, $u_n\leq c(N,q) W_{F_n}$
with  constant independent of $n$. As $\set{F_n}$
increases, $W_{F_n}\uar W_F$ everywhere in $D$. By \req{lim-Lqloc} and
\rlemma{uniformbound}, we can extract a \sseq of $\set{u_n}$ which
converges to $u$ a.e. in $D$. Hence $u\leq c W_F$. \eproof
\bcom
\bdef{gsmodlim} \normalfont{ Let $\set{\mu_n}$ be an increasing \seq
of bounded positive measures in $W^{-2,q}_+(\RN)$.
Denote
 \begin{align}\label{domainsol1}
D:=\{x\in \RN:\, \exists\, Q_x=\text{\qfine open \ngh of $x$,}&\\
\lim\mu_n(Q_x)<\infty&\}.\notag
\end{align}
Let $\mu_0;=\lim\mu_n$; the $\C2q$ {\em modified limit} of $\set{\mu_n}$, denoted by
$$\mu:=\lim^q_{n\tin}\mu_n$$
is defined as follows:
$$\mu=\mu_0 \;\text{ in $D$},\q \mu(E)=\infty \;\text{ if }\; E\cap D^c\neq \ems.$$
We also denote $u_\mu:=\lim u_{\mu_n}$. $D$ will be called {\em the regular set} of $\mu$ and $u_\mu$. }
\es
With the above notation we have,
\bth{gsmodsol} 
$\mu\in\GTM\qloc(D)$ and $u_\mu$ is a \gsmod $\C2q$-strong solution
of \req{eqmu} in $D$. \es \bproof

\eproof


\begin{equation}\label{domainsol2}\BAL
x\in D \;\iff\; &\exists\, Q_x=\text{\qfine open \ngh of $x$ \sth}\\
&u_\mu\in L^q\qloc(Q_x).
\EAL\end{equation}
 \begin{align}\label{domainsol2}
D:=\{x\in \RN:\, \exists\, Q_x=\text{\qfine open \ngh of $x$,}&\\
\lim\int_{Q_x}u_{\mu_n}^q\,dx<\infty&\}.\notag
\end{align}
\end{comment}
\bth{C(F)=0} {\rm (i)} If $F$ is a Borel set \sth $\C2q(F)=0$ then
the only non-negative $\C2q$-strong subsolution of \req{eq} in $F^c$
is the trivial solution.\1 {\rm (ii)}  If $F$ is a \qfine closed set
then $V_F$ is a \gsmod $\C2q$-strong solution in $F^c$.\1 {\rm
(iii)} Let $F$ be a \qfine closed set. If $v$
is a $\C2q$-strong solution in $D:=\RN\sms F$ then $v\leq V_F$. \es \bproof (i)
By definition,  a $\C2q$-strong solution $u$ in $D=\RN\sms F$ is the
limit of classical solutions in open sets containing $D$. In the
case that $\C2q(F)=0$, any such classical solution is the zero
solution. Hence $u=0$.\1 (ii) This is a \cons of Theorem
\ref{t:V=limU}.\1 (iii) By \rth{fine-strong} :
$$v\leq c(N,q)W_{F}\leq c'(N,q) V_{F}\; \text{ a.e. in }\;\RN\sms  F.$$
In addition, for every $\ga\geq 1$,
$$\sup\set{u: u \text { $\C2q$-strong solution in $F^c$, }u\leq \ga V_F}=V_F.$$
Hence $v\leq V_F$.
\eproof

\bth{CUcomp} Let $D$ be a \qfine open set
and let $\set{v_k}$ be a \seq of non-negative $\C2q$-strong
solutions of \req{eq} in $D$ converging a.e. in $D$. Then $v:=\lim
v_k$ is a $\C2q$-strong solution in $D$.
\es
\bproof By
\rlemma{uniformbound} there exists an increasing  \seq of compact
sets $\set{E'_n}$ \sth $\cup E'_n\sbs D$ and $\C2q(D\sms  \cup E'_n)= 0$
and $\set{v_k}$ is uniformly bounded in $L^q(E'_n)$ for
every $n$. Since  $\set{v_k}$ converges a.e. it follows that it
converges in $L^1(E'_n)$ for every $n$.  By \rth{C(F)=0} (iii)
$V_{D^c}$  dominates $\set{v_k}$.
By the dominated convergence theorem, $v_k\to v$ in the topology
$\tau^q\qloc(D)$. By assumption, for each $k$, $v_k$ is a
$\C2q$-strong solution. This means that there exists a decreasing
\seq of open sets $\set{\Gw_{m,k}}_{m=1}^\infty$, \sth
\begin{equation}\label{D-Gwmk}
 D\sbs \Gw_{m,k},\q \lim_{m\tin}\C2q(\Gw_{m,k}\sms D)=0
\end{equation}
and, for each $m$,  there exists a positive solution $u_{m,k}\in L^q\loc(\Gw_{m,k})$
of the equation $-\Gd u+u^q=0$ in $\Gw_{m,k}$ \sth
$$u_{m,k}\to v_k \txin{L^q\qloc(D)}.$$
By \rlemma{metrictop} the space $L^q\qloc(D)$ with the topology $\tau^q\qloc(D)$ is a metric space.
We denote a metric for this topology by $d\qloc$. For each $k$ let $m_k$ be sufficiently large so that
$$d\qloc(v_k, u_{m_k,k})<2^{-k} \txt{and} \C2q(\Gw_{m_k,k}\sms D)<2^{-k}.$$
Denote $v'_k=u_{m_k,k}$ and $\Gw'_k=\cap_{j=1}^{k}\Gw_{m_j,j}$. Then $\set{(v'_k,\Gw'_k)}$
is a determining \seq for $v$ in $D$.
\eproof

\bth{C(An)to0} Let $F$ be a \qfine closed set and let $\set{A_n}$ be
a  \seq of \qfine closed subsets of $F$. For each $n$, let $v_n$ be
a $\C2q$-strong solution in $D_n:=\RN\sms A_n$.

If $\C2q(A_n)\to 0$ then $v_n\to 0$ a.e. in $\RN\sms F$.

In particular, if  $\sum \C2q(A_n)<\infty$ and $v^*_n$ denotes the extension of $v_n$ to $\RN$ \sth $v^*_n=\infty$ in $A_n$
then,
\begin{equation}\label{v'nto0}
v^*_n\to 0\; \text{ a.e.  in $\RN$}.    
\end{equation}
\es
\bproof By \rth{C(F)=0}(iii)
$$v_n\leq c(N,q)W_{A_n}\leq c'(N,q) V_{A_n}\; \text{ a.e. in }\;\RN\sms  A_n.$$
By \rth{VFest} $V_{A_n}\to 0$ a.e. in $\RN\sms F$. This proves the first assertion. To verify the second assertion we apply the first
to the \seq $\set{A_n}_{n=k}^\infty$ with $F$ replaced by $F^k=\cup_{n=k}^\infty A_n$. Note that $F^k$ is \qfine closed up to a set of capacity zero.
\eproof

\bth{subsolution} Suppose that $z$ is a non-negative $\C2q$-strong
subsolution of \req{eq} in a \quasi open set $D$. Then there exists
a $\C2q$-strong solution dominating it.
\es
\bproof Let
$\set{(z_n,\Gw_n)}$ be a determining sequence for $z$. Since
$(z_n)_+$ is also a subsolution we may
 assume that $z_n\geq 0$. Let $Z_n$ be the smallest solution in $\Gw_n$ which dominates $\max(z_1,\cdots, z_n)$. Then
 $Z_n\leq Z_{n+1}$ in $\Gw_{n+1}$. Furthermore, by \rth{C(F)=0}(iii) $Z_n\leq V_{D^c}$ in $D$. Therefore, by \rth{CUcomp},
 $Z=\lim Z_n$ is a $\C2q$-strong solution in $D$.
\eproof
\bth{C2q-uniq}
Let $\Gw$ be a $\C2q$-quasi open set. Suppose that there exists a \seq of open sets $\set{G_n}$  \sth
\begin{equation}\label{un-1}\BAL
&{\rm (a)} \qq \C2q(G_n\Gd \Gw)\to 0,\\
&{\rm (b)}\qq \C2q(\prt G_n\sms \prt_q\wtl G_n)\to 0.
\EAL\end{equation}
If $v$ is a \prtql solution of \req{eq} in $\Gw$
  then $v=V_{\wtl D}$ in $\Gw$,  $D:=\RN\sms \wtl \Gw$. Thus $V_{\wtl D}$ is the unique \prtql solution in $\Gw$.
\es
\Remark Every $\C2q$-quasi open set $\Gw$ is $\C2q$-equivalent to the intersection of a \seq of open sets $\set{O_n}$
\sth $\C2q(O_n\sms \Gw)\to 0$. However, in the statement of the theorem, we do not require that $G_n$ contain $\Gw$.
Instead we require \req{un-1} (b).

The proof of the theorem is based on several lemmas. The first collects several useful formulas:
\blemma{bnd-rel}
Let $A,E_1,E_2$ be sets in $\RN$. Then the following relations hold:
\begin{equation}\label{bnd-rel}\BAL
&\mathrm{(i)} &&\prt_q A^c = \prt_q A, \\
&\mathrm{(ii)}&&\prt_q (E_1\cup E_2)\sbs \prt_q E_1\bigcup \prt_q E_2,\\
&\mathrm{(iii)} && \prt_q (E_1\cap E_2)\sbs \prt_q E_1\bigcup \prt_q E_2,\\
&\mathrm{(iv)}&&  \prt_q E_1\sbs \prt_q E_2 \bigcup \prt_q (E_2\sms E_1)\bigcup \prt_q (E_1\sms E_2),\\
&\mathrm{(v)}&& \prt_q E_1\Gd \prt_q E_2\sbs \prt_q (E_2\sms E_1)\bigcup \prt_q (E_1\sms E_2),\\
&\mathrm{(vi)}&& \prt_q A\sbs \prt A,\q \prt_q\wtl A\sbs \prt_q A.
\EAL\end{equation}
\es
\bproof (i),(ii) and (vi) follow immediately from the definition of boundary. (iii) follows from (i), (ii) and the relation
$$(E_1\cap E_2)^c=(E_1^c\cup E_2^c).$$
By (ii), 
$$\prt_q E_1\sbs \prt_q (E_1\cap E_2)\bigcup \prt_q (E_1\sms E_2).$$
By (i) and (iii), the relation,
$$E_1\cap E_2=E_2\cap(E_2\sms E_1)^c, $$
implies that
$$\prt_q(E_1\cap E_2)\sbs \prt_q E_2\bigcup \prt_q(E_2\sms E_1).$$
These relations imply (iv) which in turn implies (v).
\eproof
\note{Notation} Let $\set{A_n}$ and $ \set {B_n} $ be two sequences of sets.\1
(a) The notation $A_n\limsbs B_n$ means that $\C2q(A_n\sms B_n)\to 0$.\1
(b) The notation $A_n\limsim B_n$ means that $\C2q(A_n\Gd B_n)\to 0$.

\blemma{G-rel} Under the assumptions of the theorem,
\begin{equation}\label{G-rel}
 \prt_q G_n \limsim \prt_q \wtl G_n \limsim \prt G_n,\q \wtl G_n\limsim \ovl G_n,
\end{equation}

\begin{equation}\label{GW-rel}
 \prt_q G_n\limsim \prt_q\Gw
\end{equation}
and
\begin{equation}\label{GsmsW}
 \C2q(\ovl G_n\Gd \wtl\Gw)\to 0.
\end{equation}
In addition
\begin{equation}\label{prtqW=prtq(tlW)}
   \prt_q \Gw \qsim \prt_q \wtl\Gw.
\end{equation}
\es
\bproof By \req{un-1}(b) and \rlemma{bnd-rel} (vi) we have
$$\prt_q G_n\sbs \prt G_n\limsbs\prt_q\wtl G_n\sbs\prt_q G_n.$$
This proves \req{G-rel}.

Condition \req{un-1} (a) implies that
\begin{equation}\label{GW-bnd}
\C2q( \prt_q (G_n\sms \Gw))\to 0,\q \C2q( \prt_q (\Gw\sms G_n))\to 0.
\end{equation}
This fact and \rlemma{bnd-rel} (v) imply \req{GW-rel}.

Next observe that,
$$\BAL
\wtl G_n\sms \wtl\Gw\sbs (G_n\sms \Gw) \cup (\prt_q G_n\sms \prt_q\Gw),\q
\wtl\Gw\sms\wtl G_n\sbs (\Gw\sms G_n)\cup (\prt_q\Gw\sms\prt_q G_n).
\EAL$$
Therefore \req{un-1}(a) and \req{GW-rel} imply
\begin{equation}\label{GsmsW-2}
 \C2q(\wtl G_n\Gd \wtl\Gw)\to 0.
\end{equation}
This fact and  \req{G-rel} imply \req{GsmsW}.

By \req{bnd-rel} and  \req{GsmsW-2},
\begin{equation}\label{GsmsW-3}
\C2q(\prt_q\wtl G_n\Gd \prt_q\wtl\Gw)\to 0.
\ee
This fact together with \req{G-rel} and \req{GW-rel} imply \req{prtqW=prtq(tlW)}.

\eproof

\blemma{cap-max} Let $G$ be an open set  and  $Q$ be
 a \quasi open set. Assume that $\C2q(\prt G\sms \prt_q G)=0$.
Let $v$ be a $\C2q$-strong solution in $G'=G\sms \wtl Q$  and let
$u$ be a (classical) solution of \req{eq} in a domain $G_0$ \sth
$\bar G\sbs G_0$. Suppose that $u,v$ are non-negative and
\begin{equation}\label{u<v}
 u\qleq v \q\text{at } F:=\prt_q G\sms Q.
\end{equation}
Then
 \begin{equation}\label{u<v+V}
u\leq v+V_{\wtl Q} \txin{G'}.
 \end{equation}
\es \bproof Let $\ge$ be a positive number. Condition \req{u<v}
means that there exists  a determining \seq  $\set{(v_n,\Gw_n)}$ for
the $\C2q$-strong solution $v$ in $G'$ and  a family of open sets
$\set{Q_{k,\ge}}$  satisfying \req{Gwn-cond}, \req{Qke-cond}  and
\req{Qke-ineq F} (with $D$ replaced by $G'$). We may and shall
assume that $\Gw_n\sbs G$, that $\set{\Gw_n}$ is decreasing and
that, for every $\ge>0$, $\set{Q_{k,\ge}}_{k=1}^\infty$ is
decreasing.

In the next part of the proof we keep $\ge$ fixed. If $K$ is a compact subset of $F\sms Q_{n,\ge}$ then \req{Qke-ineq F} implies that there exists an open
\ngh of $K$, say $O_K$, \sth
$$u-v_n\leq \ge \txt{in} O_K\cap\Gw_n.$$
Therefore there exists an increasing \seq of compact sets
$\set{K_{n,\ge}}$ and a \seq of open sets $\set{O_{n,\ge}}$ \sth
\begin{align}\label{Kn+One}
  &K_{n,\ge}\sbs F\sms Q_{n,\ge}, \q \C2q(F\sms K_{n,\ge})\to 0,\\ &K_{n,\ge}\sbs O_{n,\ge},\q u-v_n\leq \ge
  \txt{in  $O_{n,\ge}\cap \Gw_n$.}\label{u-vn<e}
\end{align}
Let $\set{O'_{n,\ge}}$ be a decreasing family of open sets \sth
\begin{equation}\label{O'ne}
\BAL &(\prt G\sms \prt_q G)\cup (F\sms K_{n,\ge}) \cup \wtl Q\sbs O'_{n,\ge},\\
 &\C2q(O'_{n,\ge})\to \C2q(Q),\q \cap_{n=1}^\infty O'_{n,\ge}\qsim \wtl Q.
 \EAL\ee
Then $E_{n,\ge}:=O_{n,\ge}\cup O'_{n,\ge}$ is an open \ngh of $\prt G$ and
$$G\sms E_{n,\ge}\sbs G'\sbs \Gw_n.$$
\Consy there exist smooth domains $\Gw_{n,\ge}$ \sth
$$ \{x\in \Gw_n:\, \dist (x,\bdw_n)\geq 2^{-n}\}\sbs \Gw_{n,\ge}\sbs \ovl\Gw_{n,\ge}\sbs \Gw_n,\q
\prt \Gw_{n,\ge}\sbs E_{n,\ge}.$$

The function $w_{n,\ge}:=(u-v_n-\ge)_+$ is a classical subsolution in $\Gw_n$ and it vanishes in
$\Gw_{n}\cap O_{n,\ge}$.
Put $S_{n,\ge}=\prt\Gw_{n,\ge}\sms O_{n,\ge}$ and
$$z_{n,\ge}:=\begin{cases}w_{n,\ge}&\text{in }\bar\Gw_{n,\ge}\sms S_{n,\ge} \\0 &\text{in }
\RN\sms \bar\Gw_{n,\ge}.\end{cases}$$
Then $z_{n,\ge}$ is a (classical) subsolution in $\RN\sms S_{n,\ge}$.
Since $v_n\to v$ in $L^q\qloc(G')$, it follows that there exists a \sseq (still denoted $\set{v_n}$) \sth
 $v_n\to v$ a.e. in $G'$. Therefore $\set{z_{n,\ge}}$ converges a.e. in $D:=\RN\sms \wtl Q$ to the function
 $$z_\ge:=\begin{cases}(u-v-\ge)_+&\text{in }G'\\0 &\text{in } \RN\sms  G.\end{cases}$$
In addition
$$\sup_{\RN\sms S_{n,\ge}}z_{n,\ge}\leq \sup_{\bar G}u<\infty.$$
Note that  $D\sbs \RN\sms S_{n,\ge}$ for all $n$. Therefore, by the dominated convergence theorem,
 $z_{n,\ge}\to z_\ge$ in $L^q\qloc(D)$; \consy $z_\ge$ is a
$\C2q$-strong subsolution
in $D$. In fact $\set{(z_{n,\ge},\RN\sms S_{n,\ge})}$ is a determining \seq for $z_\ge$ in $D$. 

By \rth{subsolution}, there exists a $\C2q$-strong solution $Z_\ge$
in $D$ \sth $z_\ge\leq Z_\ge$. By \rth{C(F)=0}~(iii), $Z_\ge\leq
V_{\wtl Q}$ in $D$. Thus $z\leq V_{\wtl Q}$ and so $u-v-\ge\leq
V_{\wtl Q}$ in $G'$. Letting $\ge \to 0$ we obtain \req{u<v+V}.
\eproof
\bcom
Since $v_n\leq \ge$ in $O_{n,\ge}\cap \Gw_n$,
it follows that
$$z_{n,\ge}\leq \ge+ V\ind{A_{n,\ge}}\; \txt{ in $\Gw_{n,\ge}^c$,}$$
where
$$A_{n,\ge}=\prt\Gw_{n,\ge}\sms O_{n,\ge}\; \text { and }\;\C2q(A_{n,\ge}\sms Q)\to 0.$$
 Thus
$$z_\ge\leq \ge +V_{\wtl Q}\;\text{ in $\RN\sms G'$, }\; z_\ge=v\;\text{ in $G'$.}$$

Let $\set{(w_n,\Gw_n)}$ be a determining \seq for the $\C2q$-strong
solution $w$ in $G'$ (see \rdef{finesolution}). By definition,
$G'\sbs \Gw_n$ and we may assume that $\Gw_n\sbs G$. As $G'$ is
\qfine open, there exists a \seq of open sets $\set{G'_n}$ \sth
$G'\sbs G'_n$ and $\C2q(G'_n\sms G')\to 0$. Since we can replace
$\Gw_n$ by $\Gw_n\cap G'_n$ we may assume that $\C2q(\Gw_n\sms
G')\to 0$. As $\wtl Q\cap G'=\ems$, it follows that
$$ \C2q( \wtl Q\cap\Gw_n)\to 0.$$

Further, let $w'_n$ be the smallest solution in $\Gw_n$ dominating
the subsolution $\max(u,w_n)$. Then $\set{(w'_n,\Gw_n)}$ is a
determining \seq for $w'=$ the smallest $\C2q$-strong solution in
$G'$ \sth
$$\max(u,w)\leq w'\leq u+w.$$
Then $w'=w$ in $N_\Gg\cap G'$ and $u\leq w'$ in $G'$.

Condition \req{u<w} implies that for every compact set $K\sbs \prt G\sms Q$ there exists a positive number $\gd\ind{K}$
\sth
Let $\set{\Gw'_n}$ be a \seq of smooth domains \sth
$$\set{x\in \Gw_n: \dist (x,\bdw_n)\geq 2^{-n}}\sbs \Gw'_n\sbs \Gw_n.$$

Let $z_n$ be the largest subsolution of \req{eq} in $\RN$ 
dominated by $(u-w_n)_+$ in $\Gw'_n$.
Then $z_n=(u-w_n)_+$
in $\Gw'_n$ and $z_n$ is the solution of \req{eq} in $\RN\sms \Gw'_n$ \sth $z_n=(u-w_n)_+$ on $\bdw'_n$.
Observe that $u\leq w_n$

Define $z$
in $E:=\RN\sms \wtl Q$ as follows:
$$z:=\begin{cases}(u-w)_+&\text{in }G'\\0 &\text{in } E\sms G'.\end{cases}$$
Then
$$\int_{\RN} (-z_n\Gd \gf^2 + z_n^q\gf^2)dx\leq 0 \forevery \gf\in C_c^\infty(\RN):\; \gf\geq 0.$$

We claim that $z_n\to z$ a.e. in $E$.

Since  $z_n\to z$ in $L^q\qloc(E)$ this implies that $z$ is a
$\C2q$-strong subsolution in $E$. For this observe that $z_n\to 0$
a.e. outside $G'$ and is dominated by $u$. In addition $z_n\to z$ in
$L^q\qloc(G')$ because $w_n\to w$ there. Thus $z_n\to z$ a.e. and is
dominated by $u$.
\end{comment}
\bcor{cap-max} Let $G,Q,G'$ and $u,v$ be as in the statement of the lemma.
If
\begin{equation}\label{blowup-prt G'}
 \qlim_{x\to \prt_q G'}v(x)=\infty \q\text{\qae at}\; \prt_q G'
\ee
then \req{u<v+V} holds.
\es
\bproof Since $u$ is bounded in $G$, \req{blowup-prt G'} implies \req{u<v}.
Therefore the previous lemma implies \req{u<v+V}.
\eproof

\bcom
Let
$$A_w:=\set{y\in  \prt G\sms \wtl Q:\,\liminf_{G' \ni x\to y}w(x)<\infty}.$$
Since $Q$ is \quasi open and, by assumption, $\C2q(A_w)=0$ it follows that, for every $\ge>0$, there exists an open set
 $Q'_\ge$ \sth
$$Q\cap A_w\sbs Q'_\ge,\qq \C2q(Q'_\ge\sms Q)<\ge.$$
Given $\gd>0$ put
$$G_\gd:=\set{x\in G:\dist(x,\prt G)>\gd},\qq G_{\gd,\ge}:=G_\gd\sms \wtl Q'_\ge.$$
For every  $\ge>0$ there exists $\gd_\ge>0$ \sth $u\leq w$ on $\prt G_\gd\sms Q'_\ge$ for every $\gd<\gd_\ge$.
Therefore, by \rlemma{cap-max},
$$u\leq w+V_{\wtl{Q'}_\ge} \txin{G_{\gd,\ge}}.$$
As $\gd\dar 0$, $G_{\gd,\ge}\dar G\sms \wtl Q'_\ge$; \consy
$$u\leq w+V_{\wtl {Q'}_\ge} \txin{G\sms \wtl Q'_\ge}.$$
 Further,
$$V_{\wtl {Q'}_\ge}\leq V_{\wtl {Q'}}+ V_{\wtl {Q'}_\ge\sms \wtl {Q'}}$$
and $\wtl {Q'}_\ge\sms \wtl {Q'}\sbs \wtl{ {Q'}_\ge\sms Q'}$. Since $\lim_{\ge\to 0}\C2q(\wtl{ {Q'}_\ge\sms Q'})=0$
it follows, by \rth{C(F)=0}, that $\lim_{\ge\to 0}V_{\wtl {Q'}_\ge\sms \wtl {Q'}}=0$. This implies \req{u<w+V}.
\end{comment}

\note{Proof of \rth{C2q-uniq}} Let $K$ be a compact subset of $D=\RN\sms \wtl\Gw$ and let $\mu\in W^{-2,q}(\RN)$ be a non-negative
measure supported in $K$. We  prove that
\begin{equation}\label{C2q-0}
   u_\mu\leq v \txin{\Gw}
\end{equation}

As $\mu(\wtl\Gw)=0$, \req{GsmsW} implies
\begin{equation}\label{C2q-1}
\mu'_n:=\mu\chr{\ovl{G}_n}=\mu\chr{\ovl G_n\sms\wtl\Gw} \to 0.
\end{equation}
If $\nu$ is a bounded measure \sth $\nu(\ovl G_n)=0$ and $O_{n,k}$ is a \seq of open \nghs of $\ovl G_n$ \sth
$\cap_kO_{n,k}=\ovl G_n$ then
$$\nu(O_{n,k}\sms \ovl G_n)\to 0 \txt{as $k\tin$.}$$
Applying this observation to $\nu=\mu-\mu'_n$ and using \req{C2q-1} we conclude that,
 for every $n\in \BBN$, there exists a non-negative measure
$\mu_n$ \sth
\begin{equation}\label{C2q-2}
 \mu_n\leq\mu,\q \supp\mu_n\cap\ovl G_n=\ems, \q (\mu-\mu_n)(\RN)\to 0.
\end{equation}
As $K_n:=\supp \mu_n$ is a compact set disjoint from $\ovl G_n$ it follows that $u_{\mu_n}$ is a bounded solution
of \req{eq} in a \ngh of $\ovl G_n$.

Let $Q_n:=G_n\sms\wtl \Gw$. By  \req{GsmsW} $\C2q(Q_n)\to 0$ and therefore
\begin{equation}\label{C2q-3}
   \C2q(\wtl Q_n)\to 0.
\end{equation}
Applying \rcor{cap-max} to $G_n,Q_n$ with $u=u_{\mu_n}$ and $w=v$ we obtain
\begin{equation}\label{C2q-4}
  u_{\mu_n}\leq v+V_{\wtl Q_n} \txin{G_n\sms\wtl Q_n}.
\end{equation}
By \req{C2q-2} $u_{\mu_n}\to u_\mu$ and, by \req{C2q-3}, $V_{\wtl Q_n}\to 0$. Therefore, in view of \req{un-1}~(a),
\begin{equation}\label{C2q-5}
  u_{\mu}\leq v \qq \qae \text{ in $\Gw$}.
\end{equation}
This holds for every non-negative measure $\mu\in W^{-2,q}(\RN)$ supported in a compact subset of $D$.
Therefore 
\begin{equation}\label{C2q-6}
  V_{\wtl D}=V_D\leq v \qq \qae \text{ in $\Gw$}.
\end{equation}
On the other hand, by \rth{C(F)=0}~(iii), $v\leq V_{\RN\sms \Gw}$. But \req{prtqW=prtq(tlW)} implies
that
$$\RN\sms \Gw=D\cup\prt_q\Gw\qsim D\cup\prt_q\wtl\Gw.$$
As $\prt_q D= \prt_q\wtl \Gw$ it follows that
$\RN\sms \Gw\qsim \wtl D.$
Thus $V_{\wtl D}=V_{\RN\sms \Gw}$ and finally $v=V_{\wtl D}$.

\qed

\note{Example} Let $\set{x^m}$ be a \seq of distinct points in $B_1(0)$. Let $\set{r_n}$ be a decreasing \seq of positive
numbers \sth \req{cor-un2} holds and $\ovl B_{r_n}(x^n)\sbs B_1(0)$. Put
$$\Gw_n=B_1(0)\sms \cup_1^n \ovl B_{r_k}(x^k).$$
Then there exists a unique large solution in
$$\Gw:=\cap_1^\infty \Gw_n=B_1(0)\sms \cup_1^\infty \ovl B_{r_k}(x^k).$$

\appendix
\section{On the space $W^{2,q'}_{0,\infty}$}\label{w0q}
We establish some features of the space $W^{2,q'}_{0,\infty}$ which show that it is sufficiently rich in order to serve
as a space of test functions in \qfine open sets. These are used mainly in Section \ref{fineq}.

\blemma{qsmoothing} Suppose that  $D$ a $\C2q$ finely open set and $K$ is a bounded \qfine closed subset of $D$.
Then, for every $a>0$, there exists $\gf_a\in W^{2,q'}(\RN)$ \sth:
\begin{equation}\label{qchar}\BAL
 &(i)\q 0\leq\gf_a\leq 1,\q(ii)\q \qsupp \gf_a \Sbs D,\\&(iii)\q\C2q(\set{x\in K:\,\gf_a(x)<1})<a.
\EAL\end{equation}

\es
\bproof Let $0<\ge(1+2^{q'})<a$. Let $K'$ be a compact set  and  $D'$  be an open set \sth,
$$ K'\sbs K,\q D\sbs D',\q \C2q(K\sms K')<\ge,\q \C2q(\widetilde{D'\sms D})<\ge.$$
 Let $\gf$ be a smooth function with compact support in $D'$
\sth $0\leq \gf\leq 1$ and $\gf=1$ on a \ngh of $K'$. Let $\set{A_n}$ be a decreasing \seq of open \nghs
of $\widetilde{D'\sms D}$ \sth
$$\C2q(\tl A_n)\to\C2q(\widetilde{D'\sms D}).$$
 Further, let $\set{\eta_n}$ be a \seq of functions in $\W2q(\RN)$ \sth
$$0\leq \eta_n,\q  \eta_n\geq 1 \text{ \qae in }{\tl A_n},\q   \norm{\eta_n}^{q'}\ind{\W2q(\RN)}=\C2q(\tl A_n).$$
(See \cite[Thm.2.3.10]{AH} for the existence of such functions.)

Let $\ga\in (0,1)$ and put $E_n=\set{x\in D:\, \eta_n(x)\geq 1-\ga}$. Then
$$(1-\ga))^{-q'}\norm{\eta_n}\ind{W^{2,q'}(\RN)} \geq \C2q(E_n)$$
so that
$$\limsup \C2q(E_n)\leq \C2q(\tl A_n) /(1-\ga)^{q'}<\ge /(1-\ga)^{q'}.$$
Let $h$ be a monotone, smooth cutoff function \sth
\[h(t)= \begin{cases} 0 &\text{if }t<\ga/4\\ h(t)=t &\text{if }t>\ga/2.\end{cases}\]
 Then
 $\gf_n:=h\circ(\gf-\eta_n)\in \W2q(\RN)$ and
$$\gf_n\geq \ga \text{ on } K'_n:=K'\sms E_n,\q  \gf_n=0 \text{ \qae in } A_n.$$
Thus, choosing $\ga =1/2$,
\begin{equation}\label{char1}\begin{cases}
 &\gf_n/\ga\geq 1 \text{ on } K'_n,\q \qsupp\gf_n\sbs (\supp\gf)\sms A_n \Sbs D, \\ &\limsup\C2q(K\sms K'_n)<\ge(1+(1-\ga)^{-q'}).
\end{cases}\end{equation}
By applying (to $\gf_n/\ga$)  another smooth
cutoff function which approximates $\min(\cdot,1)$, we obtain a \seq  of functions which, for $n$ sufficiently large,
satisfy the statement of the lemma.
\eproof
\bcor{partition1} Let $D$ be a bounded \qfine open set and let $\set{E_n}$ be a $q$-exhaustion of $D$
(see \rlemma{qexhaust}). Then there exists a \seq $\set{\vgf_n}$ in $\W2q(\RN)$ \sth:
\begin{equation}\label{partition1}\BAL
   &(i)\;0\leq \vgf_n\leq 1,&&(ii)\; \qsupp\vgf_n\Sbs E_{n+1},\\
   &(iii)\;\sum_{n=1}^\infty\C2q(E_n\sms [\vgf_n=1])<\infty,&& (iv)\; \set{\vgf_n} \text{ is non-decreasing.}
\EAL\end{equation}
In particular $\vgf_n\uar 1$ \qae in $D$.
\es
\bproof We construct $\phi_n$ as in \rlemma{qsmoothing} with $K$ and $D$ replaced by $E_n$ and
$\mathrm{int}_q E_{n+1}$,
$a=2^{-n}$ and $\ga=1/2$.
Then we put $\tl\vgf_n:=\sum_1^n\vgf_m$ and finally  apply to  $2\tl\vgf_n$ a smooth
cutoff function which approximates $\min(\cdot,1)$.
\eproof
\blemma{metrictop} Let $D$ be a \qfine open set and let
$\tau\qloc(D)$ be the topology in $L^q\qloc(D)$ defined by the
family of seminorms \req{seminorms}. Then $\tau\qloc(D)$ is a metric
topology.
\es
\bproof It is sufficient to show that the space is separable. For each fixed $\gf\in\w02q(D)$, the space
$L^q(D;\gf)$ is separable. Let $\set{\vgf_n}$ be as in \rcor{partition1}. Then, for every $f\in L^q\qloc(D)$,
$$\int_D f\psi(1-\vgf_m)dx\to 0 \forevery \psi\in \w02q(D).$$
Therefore,
if $\set{h_{k,m}}_{k=1}^\infty$ is a dense set in $L^q(D;\vgf_m)$ then
$$\set{h_{k,m}:\, k,m\in \BBN}$$
is a dense set in $L^q\qloc(D)$.
\eproof

\blemma{uniformbound}  Assume that $F$ is a \qfine closed set and $F\sbs B_{R/2}(0)$.
Put $D=B_R(0)\sms F$. Let $\set{E_n}$ be a
$q$-exhaustion of $D$. Then there exists a $q$-exhaustion $\set{E'_n}$ \sth
\begin{equation}\label{ubound}
E'_n\sbs E_n, \q \C2q(E_n\sms E'_n)\to 0,
\end{equation}
for which the following statement holds:

The set of non-negative very weak subsolutions of \req{eq} in $D$
is uniformly bounded in $L^{q}(E'_n)$ for every $n\in \BBN$. \es
\bproof Let
$\set{\vgf_n}$ be as in \rcor{partition1} and let $A_{n,k}$ be an
open \ngh of $E_n\sms [\vgf_n=1]$ \sth
$$\C2q(A_{n,k})\leq (1+2^{-k})\C2q(E_n\sms [\vgf_n=1]),\q \tl A_{n,k+1}\sbs A_{n,k} \forevery k\geq n,\,n\in \BBN.$$
Put $E'_n=E_n\sms \cup_{k=n}^\infty A_{n,k}$. Then $\set{E'_n}$ is a $q$-exhaustion of $D$ and \req{ubound} holds.
Furthermore,
$\vgf_n=1$ on $E'_n$. Hence,
by \rth{int<W2qnorm}, every non-negative very weak subsolution $u$ of \req{eq} in $D$ satisfies
$$\int_{E'_n}u^q\,dx \leq  \int_D u^q\vgf_n^{2q'}dx\leq c(q,N)\norm{D^2\vgf_n}\ind{L^{q'}(D)}^{q'}.$$
\eproof

\blemma{Lpq} Let $K$ be a bounded \qfine closed subset of $D$. Then, for every $\ge>0$,
there exists a compact set $K_\ge\sbs K$ \sth
$$ \C2q(K\sms K_\ge)<\ge,\q f\in L^q(K_\ge) \forevery f\in L^q\qloc(D)
$$
\es
\bproof This is an immediate consequence of \rlemma{qsmoothing}.
\eproof
\section{Open problems} There are many interesting problems related to possible {\em extensions} of the theory of
solutions in finely open sets, presented in Section \ref{fineq}. We do not describe here problems of this nature, but only
problems directly related to results presented in the present paper.

 In order to formulate the first problem, it is convenient to introduce an additional definition.
\bdef{Cweak} {\rm Let $u$ be a non-negative measurable function in a \qfine open set $D$ and let
$\mu\in\GTM\qloc(D)$ be a non-negative measure.
We say that $u$ is a
{\em $\C2q$-weak solution} of \req{eqmu} in $\Gw$ if $u$ satisfies \req{weaksol} and \req{weakeq}.}
\es
\noindent{\bf Problem I.}\hskip 2mm We know that if $u$ is a $\C2q$-strong solution then it is also a $\C2q$-weak
solution, (see \rth{fine-strong}). Does the opposite implication hold: is it true that every
$\C2q$-weak solution of \req{eqmu} is a $\C2q$-strong solution?\2
\noindent{\bf Problem II.}\hskip 2mm This problem is related to
\rth{CUcomp}. The question is if the following related assertion is  valid:\2
{\em  Let $D$ be a \qfine open set
and let $\set{v_k}$ be a \seq of non-negative $\C2q$-strong
solutions of \req{eq} in $D$. Then there exists a \sseq $\set{v_{k_j}}$, converging in $L^q\qloc(D)$.}\1
We observe that if such a \sseq exists then one can extract a further \sseq which converges a.e. in $D$
and, by \rth{CUcomp},  its limit is a $\C2q$-strong solution in $D$.\2
\indent The next problem is related to the uniqueness result \rlemma{almostlarge 1}.
It is known that in the subcritical case condition \req{un-open4} is necessary in order to guarantee
uniqueness of large solutions. (In the subcritical case the notion of 'large solution' and '\prtql solution' coincide.)
The situation is essentially different \wrto \prtql solutions in the supercritical case.
In fact it is likely
that condition \req{un-open4} is not necessary in this case.\1
\noindent{\bf Problem III.}\hskip 2mm {\em Let $\Gw$ be a bounded open set and put $F:=\RN\sms \Gw$. We know that $V_F$ is an almost large
solution and, a-fortiori, a \prtql solution in $\Gw$.
Question: Is $V_F$ the unique \prtql solution in $\Gw$?}

\end{document}